\documentclass[a4paper,12pt,twoside]{article}
\usepackage{graphicx}
\usepackage[section]{placeins}
\usepackage[utf8]{inputenc}
\usepackage{amssymb}
\usepackage{amsmath}
\usepackage{amsthm}
\usepackage[linesnumbered,ruled,vlined]{algorithm2e} 
\usepackage{mathtools}
\usepackage[a4paper,text={150mm,255mm},centering,headsep=5mm,footskip=10mm]{geometry}
\usepackage[onehalfspacing]{setspace}
\usepackage[ngerman, english]{babel}
\usepackage{bbm}
\usepackage{grffile} 
\usepackage{verbatim}
\usepackage{caption}
\usepackage{subcaption}
\usepackage{xcolor}
\usepackage{gensymb}
\usepackage{tcolorbox}
\usepackage{enumitem}
\usepackage{colortbl}
\usepackage{fancyhdr}
\usepackage{url}
\usepackage[noadjust]{cite}
\usepackage{titlesec}
\usepackage{graphicx}
\usepackage{colortbl}
\usepackage[]{algorithm2e}
\usepackage{hyperref}
\usepackage{booktabs}
\usepackage[arrow, matrix, curve]{xy}
\usepackage{authblk}
\usepackage[font=small]{caption}
\usepackage{stackengine} 
\usepackage{tikz}
\usetikzlibrary{arrows,arrows.meta}
\usetikzlibrary{calc,trees,positioning,arrows,chains,shapes.geometric,%
    decorations.pathreplacing,decorations.pathmorphing,shapes,%
    matrix,shapes.symbols}

\tikzset{
>=stealth',
  punktchain/.style={
    rectangle, 
    rounded corners, 
    draw=black, very thick,
    text width=10em, 
    minimum height=3em, 
    text centered, 
    on chain},
  line/.style={draw, thick, -},
  element/.style={
    tape,
    top color=white,
    bottom color=blue!50!black!60!,
    minimum width=8em,
    draw=blue!40!black!90, very thick,
    text width=10em, 
    minimum height=3.5em, 
    text centered, 
    on chain},
  every join/.style={-, thick,shorten >=pt},
  decoration={brace},
  tuborg/.style={decorate},
  tubnode/.style={midway, right=0pt},
}

\tcbuselibrary{breakable}

\definecolor{mycolor}{rgb}{0,0,0}
\makeatletter
\newcommand{\greybox}[1]{%
  \setbox0=\hbox{#1}%
  \setlength{\@tempdima}{\dimexpr\wd0+13pt}%
  \begin{tcolorbox}[colframe=mycolor,breakable,boxrule=0.2pt,arc=4pt,colback=white!10,
      left=6pt,right=6pt,top=6pt,bottom=6pt,boxsep=0pt,width=
\textwidth] 
    #1
  \end{tcolorbox}
}

\newtheoremstyle{boldremark}
    {\dimexpr\topsep/2\relax} 
    {\dimexpr\topsep/2\relax} 
    {}          
    {}          
    {\bfseries} 
    {.}         
    {.5em}      
    {}          

\theoremstyle{remark}

\theoremstyle{boldremark}

\newcommand*{\R}{\mathbb{R}}

\newcommand*{\fatX}{\textbf{X}}
\newcommand*{\fatx}{\textbf{x}}

\DeclareMathOperator*{\argmin}{arg\,min}

\newcommand*{\vertbar}{\rule[-1ex]{0.5pt}{2.5ex}}

\newcommand{\pk}[1]{\textcolor{red!85!black}{\textbf{PK: }#1}}
\newcommand{\nw}[1]{\textcolor{blue!85!black}{\textbf{NW: }#1}}
\newcommand{\MYhref}[3][blue]{\href{#2}{\color{#1}{#3}}}%
\numberwithin{equation}{section}
\numberwithin{theorem}{section}
\numberwithin{definition}{section}
\numberwithin{lemma}{section}
\numberwithin{corollary}{section}
\numberwithin{remark}{section}

\usepackage{geometry}
\titleformat{\chapter}[display]
   {\normalfont\huge\bfseries}{ \thechapter}{20pt}{\Huge}

\titlespacing*{\chapter}{0pt}{-40pt}{20pt}
\definecolor{RoyalRed}{RGB}{157,16,45}
\titleformat{\chapter}[display]
  {\bfseries\LARGE}
  {\huge
  \chaptertitlename\hspace{0.1ex} \thechapter}{1pc}
  {{\titlerule[0pt]}\vspace{1pc}}

\titleformat{\chapter}[hang]{\bfseries\huge}{\bfseries\thechapter}{1em}{}

\usepackage[scaled]{helvet}
\usepackage{courier}

\providecommand{\keywords}[1]
{
  \small	
  \textbf{\text{Keywords---}} #1
}

\begin{document}
\title{Memory-based reduced modelling and data-based estimation of opinion spreading}

\author[1]{Niklas Wulkow}
\author[1]{Péter Koltai}
\author[1,2]{Christof Schütte}
\affil[1]{Department of Mathematics and Computer Science, Freie Universität Berlin, Germany}
\affil[2]{Zuse Institute Berlin, Germany}
\maketitle

\begin{abstract}
  We investigate opinion dynamics based on an agent-based model, and are interested in predicting the evolution of the percentages of the entire agent population that share an opinion. Since these opinion percentages can be seen as an aggregated observation of the full system state, the individual opinions of each agent, we view this in the framework of the Mori--Zwanzig projection formalism.
  More specifically, we show how to estimate a nonlinear autoregressive model (NAR) with memory from data given by a time series of opinion percentages, and discuss its prediction capacities for various specific topologies of the agent interaction network. We demonstrate that the inclusion of memory terms significantly improves the prediction quality on examples with different network topologies.
  \end{abstract}

\setcounter{tocdepth}{3}

\section{Introduction}

Political opinion polls capture how the opinions of people within a society regarding a certain topic or their current voting preferences are distributed. Individual opinions do not have to be constant but rather are subject to change induced by impactful events or the opinions of their peers which is formalized under the term \textit{conformity} in~\cite{stangor}. There have been recent advances in simulating the process in which members of a society change their opinions,  see, e.g., \cite{banisch,klimek,misra,qian,nardini,boehme,bolzern} and the review articles \cite{anderson,haoxing,castellano,sirbu}. This is in part due to increasing computing power which enables to carry out agent-based models that simulate behaviour of members of a synthetic population, such as members of a society, on the microscale by emulating the decision-making rules. The agents are often treated as the nodes of a network while an edge between two nodes means that these agents are neighbors of each other and thus influence each other's respective opinions.

One is often not interested in modelling, or predicting, which person has which opinion but rather, as in polls, what the percentage of each opinion within the society is. There is ample interest in deriving dynamics for the evolution of these percentages. 

In this article, we will present a framework which identifies the governing equations for the dynamics of opinion percentages for different types of networks, more precisely, how the governing equations can be inferred from data on the opinion percentages. To this end, we will emulate the decision making process with a simple agent-based model (ABM) that is based on the assumption of conformity and inspired by the ABM in \cite{misra}. Introductions into agent-based modelling in general can be found in \cite{jennings} and \cite{Laubenbacher2009} and specifically into agent-based models for opinion dynamics in \cite{banischbook}.

The literature contains a variety of approaches for finding governing equations on the macro-level (here: opinion percentages) based on micro-dynamics (here: agent-based model). However, most do not deal with opinion formation or voter models but with models originating from the context of the natural sciences. There it is well-known that the aggregation process from the micro- to the macro-level typically leads to non-Markovian processes, i.e., finding the governing equations on the macro-level requires the inclusion of memory, cf.\ the Mori--Zwanzig formalism \cite{zwanzig, linfu, chorin2002}. In the context of opinion formation, this aspect is hardly discussed at all. Banisch  \cite{banisch2014} discusses the issue for agent-based models; he gives stochastic and combinatorial arguments for the appearance of memory with heterogeneous micro-structure but does not present any practical methods for finding appropriate governing equations for the macro-dynamics. Several other authors discuss the micro-macro-aggregation problem in opinion formation, e.g., via influence matrices between agents \cite{scaglione,ravazzi,abir} but ignore memory effects entirely. Others discuss memory effects, but only on the micro-level, e.g., \cite{jedre,chen} (agents have memory), \cite{moussaid} (agents gain experience), or \cite{boschi} (micro-dynamics depends on collective memory). Very few articles consider the practical methods for finding governing equations on the macro-level, e.g., by inferring them from micro-level simulation data, but memory effects are ignored, cf.~\cite{lu}. Thus, there is a significant gap between Banisch' insight that opinion aggregation introduces memory and its practical use for finding appropriate description of the resulting macro-dynamics. 

This article aims at closing this gap by (1) utilizing techniques like the Mori--Zwanzig formalism and Taken's well-known embedding theorem for showing that agent-based models for the micro-dynamics lead to memory effects on the macro-level if the interaction between the agents is heterogeneous, while doing this in a way that allows for (2) proposing practical algorithmic techniques to learn governing equations for the macro-dynamics including memory utilizing macro-observations of micro-level simulation data.

More precisely, we investigate complete and incomplete interaction networks: in complete networks every agent interacts with all others (homogeneous interaction), while in incomplete networks there are subcommunities within the society that have few links between each other (heterogeneous interaction). As we will show, in the case of a complete network, one can identify a \emph{Markovian} model for the macro-dynamics of the opinion percentages using standard well-mixedness arguments known from the mean-field approaches or population limits, e.g., for predator-prey models~\cite{berryman}. However, arguments used for that case do not hold true in cases when the network is not complete. We will show how to use information from the past (memory) via a kind of delay embedding of the dynamics to describe the evolution of opinion percentages in the general case.

The exact reason for the inclusion of memory will formally by derived in Section \ref{sec:MZ} by using the Mori--Zwanzig formalism \cite{zwanzig, linfu, chorin2002}. Inspired by problems in statistical physics, the Mori--Zwanzig formalism explains how in the case of only low-dimensional observations of a high-dimensional system being available, the evolution of these observations of the full system can be obtained by replacing the missing information of the full system by past information of these available observations. This is in light of the result of Takens \cite{takens} that states that, under fairly generic assumptions, the delay embedding of the dynamics of an observable is diffeomorphic to the dynamics of the full system.

There are various techniques for the modelling of time-discrete dynamical systems which involve the memory of the system. An intuitive approach is comprised by Higher-order Markov Models \cite{raftery,tuyen}. These models are defined by transition probabilities between discrete states where each state represents a sequence of cells of a discretization of the state space with a given length (``memory depth''). Although these models can be powerful in investigating the long-term behaviour of the process by means of Markov State Models for Markovian processes \cite{noemsm}, they yield two problems: The loss of accuracy obtained from the discretization and an exponentially increasing number of states with increasing length of the sequences and number of grid cells.

Another example is simplex projection as in \cite{sugihara90} where, using Takens' result, subsequent states of a system are predicted from relative next steps of similar patterns as its recent history. 
A younger modelling technique is Long short-term memory neural networks (LSTMs) \cite{hochreiter,shaowu} which is a subclass of recurrent neural networks and specifically designed for prediction of time series for which past information is vital. However, both these techniques provide little to no understanding of the dynamical \emph{rules} of the system: Simplex projection does not produce any model or dynamical law but rather uses a procedure similar to the nearest neighbors classification algorithm (see, e.g., \cite{DGL13}). LSTMs, as most neural networks, typically have far too many parameters to admit interpretability. An additional means for forecasting of memory-dependent dynamical systems is the well-known class of autoregressive (AR) models \cite{brockwell}, which describes the evolution of a system by a linear combination of its most recent states. Additionally, there exist variants of these AR models that are sparse \cite{davissparse,fugita} or nonlinear \cite{billings} or comprise both aspects in application to a Singular Value Decomposition of a data matrix~\cite{havok}.
As we will see, linear (Markovian) systems cannot describe the evolution of opinion percentages even in the simplest case,
but \emph{simple} polynomial terms are sufficient for fully connected networks. We shall address this point with nonlinear AR (NAR) models, as derived through the Mori--Zwanzig formalism.

In addition to the analysis of micro-macro-aggregation for opinion formation, further novelty in our work lies in the methods we propose for learning NAR models from data, to describe the evolution of opinion percentages, and their theoretical justification. 
We will show that the prediction accuracy of the NAR models for the opinion percentages increases with larger memory depths.
To this end, we will deploy methods from data-driven (sparse) system identification---as in Dynamic Mode Decomposition \cite{schmid, dmd,jovanovic} or Sparse Identification of Nonlinear Dynamics (SINDy) \cite{sindy}---to the field of opinion dynamics. More precisely, we will extend SINDy towards finding (sparse) NAR models to describe the evolution of opinion percentages. The new method is called ``Sparse Identification of Nonlinear Autoregressive Models'' (SINAR), as it is technically a natural generalization of SINDy by including non-linear memory terms.
We will demonstrate that SINAR is well-suited for our purposes in learning macroscopic opinion dynamics. A conceptually similar method has been introduced in \cite{havok} with Hankel Alternative View Of Koopman (HAVOK). It can be interpreted as a special case of SINAR.

\paragraph{Outline.} In Section \ref{sec:MZ} we start with outlining the opinion aggregation process and proceed with the derivation of NAR models for the evolution of observations through the Mori--Zwanzig formalism. Next, in Section \ref{sec:sinar}, we present the SINAR method for estimating the coefficients in these NAR models from data. Last, we demonstrate how to apply SINAR for increasing the accuracy of prediction of opinion percentages in the case of incomplete interaction networks in Section~\ref{sec:abm}.

\section{Derivation of a nonlinear autoregressive model using the Mori--Zwanzig formalism}
\label{sec:MZ}

Below, we will model the spread of opinions inside a closed society by an agent-based model. It will consist of a high number $N$ of agents who change their opinions $X_i$, $i=1,\ldots, N$, within a finite set of $M$ possible opinions over discrete time steps according to a rule that is based on the opinions of themselves and other agents. This rule will be Markovian, or memory-free, i.e., the changes of opinions are only influenced by opinions in the current time step. These dynamics will be called the \textit{microdynamics}. The state of the microdynamics at time $t$ is denoted by $X_t = [(X_t)_1,\dots,(X_t)_N]^T$. The respective state space is denoted by~$\mathbb{X}$, and has cardinality $|\mathbb{X}| = M^N$.

We will only be able to observe the percentages of opinions, i.e., the ratios of those among all agents with each of the~$M$ opinions. In this article, we are interested in identifying the dynamical rules of the evolution of the percentages of opinions, which we call the \textit{macrodynamics}. Identifying the dynamics of low-dimensional observations of a higher dimensional system is a typical setup for the Mori--Zwanzig formalism~\cite{zwanzig, chorin2002, linfu}. We will consider a general framework for this and show how it yields a nonlinear autoregressive model \cite{billings} for the macrodynamics. Later on we show how it can be applied to the specific case of the spread of opinions.

\subsection{The setting: Microdynamics and projected observations}
First we assume that the microdynamics are Markovian (memory-free) and deterministic.
We consider the dynamical system $F:\mathbb{X}\rightarrow \mathbb{X}$ that governs the microdynamics
\begin{equation}
X_{t+1} = F(X_t) \in \mathbb{X}.
\label{eq:dynSystem}
\end{equation}
Further, we denote the space of observations of the microdynamics (observables) by $\mathbb{Y} \subseteq \mathbb{R}^m$ and by $\mathcal{G} := \lbrace g : \mathbb{X}\rightarrow \mathbb{Y} \rbrace$ the set of functions that map states of the dynamical system \eqref{eq:dynSystem} to $\mathbb{Y}$. We suppose from here on that we do not have knowledge of the state of the microdynamics at any point in time but instead only have the value of the fixed observable $x = \xi(X) \in \mathbb{Y}$ which we call the the accessible, or \textit{relevant}, variables.

Additionally, we define the subspace $\mathcal{H}$ of functions in $\mathcal{G}$ that depend only on these relevant variables and map to $\mathbb{Y}$ as
$\mathcal{H} := \lbrace h \in \mathcal{G}\mid \exists \tilde{h}:\xi(\mathbb{X}) \rightarrow \mathbb{Y} :\ h = \tilde{h}\circ \xi \rbrace$.
Functions in $\mathcal{H}$ still depend on $X \in \mathbb{X}$ but the information of $\xi(X)$ is enough to evaluate them. When we write $h(x)$ for $x \in \mathbb{Y}$, we abuse notation and mean $h(\xi(x))$. An example is
\begin{equation*}
\mathbb{X} = \R^2,\quad \xi(X) = X_1+X_2,\quad h(X_1,X_2) = (X_1+X_2)^2 = \xi(X)^2.
\end{equation*}
In this case it is enough to know the value of $\xi(X)$ to evaluate $h(X)$.

The goal is now to represent the evolution of the observations $x_t = \xi(X_t)$ under the microdynamics with knowledge only about values of $x_t$ but not of the states $X_t$ of the microdynamics. As illustrated in the following diagram, instead of taking one step of the microdynamics and then evaluating $\xi$, we only have access to the observation $\xi(X)$ and want to evaluate $\xi(F(X))$ under the premise that $\xi(X) = x$.
\begin{equation}
\begin{xy}
  \xymatrix{
      X \ar[r]^{F} \ar[d]_{\xi}    &   F(X) \ar[d]^{\xi}  \\
      \xi(X) = x \ar[r]_{?}       &   \xi(F(X))  
  }
\end{xy}
\label{eq:MZdiagram2}
\end{equation}
To this end, we define a projection operator $P:\mathcal{G}\rightarrow \mathcal{H}$ that maps a function depending on $X$ to a function depending on $\xi(X)$. We additionally define its complement $Q := Id - P$. We assume from now on that the microdynamics are stationary with an $F$-invariant probability distribution $\mu$ over $\mathbb{X}$, so that when asking what $g(X)$ is, we assume that $X_t$ is distributed by~$\mu$. \footnote{A natural candidate for $P$ would be the conditional expectation with respect to $\mu$, given by
$(Pg)(x) = \mathbb{E}[g(X) \mid \xi(X) = x]$,
see Appendix~\ref{sec:appendixCondExpec}. 
Approximating the conditional expectation can be a challenging task, see \cite{gilani}. Instead, we consider the orthogonal projection onto basis functions since we are seeking models spanned by such functions with the option to control the sparsity of the model. In \cite{chorin2002}, the connection between both projections is discussed.
}
We, of course, are interested in the case $g = \xi\circ F$.
We follow \cite{linfu} until the end of Section \ref{sec:formulatingMacro} and define $P$ as the orthogonal projection onto the span of a set of linearly independent functions from $\mathcal{H}$. These functions are denoted by $\varphi_1,\dots,\varphi_L:\mathbb{Y}\rightarrow \R^{m}$ which build the columns of $\varphi = [\varphi_1,\dots,\varphi_L]$.
\begin{equation}
(Pg)(x) :=  \varphi(x) \langle\varphi,\varphi\rangle^{-1} \langle\varphi,g\rangle
\label{eq:orthProj}
\end{equation}
where $x \in \mathbb{Y}$ and the scalar product $\langle\cdot,\cdot\rangle$ is defined for matrix-valued functions $f:\mathbb{X}\rightarrow \R^{m \times a}$ and $g:\mathbb{X}\rightarrow \R^{m \times b}$ as 
\begin{equation*}
\langle f,g\rangle := \int_{\mathbb{X}} \underbrace{f(X)^T}_{\in \R^{a\times m}} \underbrace{g(X)}_{\in \R^{m \times b}} d\mu(X) \in \R^{a \times b},
\end{equation*}
which itself is matrix-valued. The term $\langle\varphi,\varphi\rangle$ is a \text{mass matrix} that ensures that $P$ is an orthogonal projection. This orthogonal projection has the property that $Pg$ is the closest function in $span(\varphi)$ to $g$ with respect to $\langle \cdot, \cdot \rangle$.

Note that if $\mathcal{H}$ is infinite-dimensional, one would need an infinite number of functions to yield that $span(\varphi) = \mathcal{H}$. In this case the projection formalism is well defined if $\mathcal{H}$ is closed. In practice, in this case for the computation that will follow one would choose a sufficiently rich finite set of functions so that $span(\varphi)$ covers those parts of $\mathcal{H}$ that are of interest.

\subsection{Mori--Zwanzig representation of the macrodynamics}
\label{sec:formulatingMacro}
We will now show how to represent the evolution of the observations over time. With the Koopman operator \cite{koopman} $\mathcal{K}$ for the system \eqref{eq:dynSystem}, defined as the operator that maps a function $g \in \mathcal{G}$ to $g\circ F \in \mathcal{G}$, we consider the Dyson formula 
\begin{equation}
\mathcal{K}^{t+1} = \sum_{k=0}^t \mathcal{K}^{t-k} P\mathcal{K} (Q\mathcal{K})^k + (Q\mathcal{K})^{k+1}.
\label{eq:dyson}
\end{equation}
The Dyson formula describes a way to iteratively split up the application of the Koopman operator to a function $g$ into parts $P\mathcal{K}g$ and $Q\mathcal{K}g$. Equation \eqref{eq:dyson} yields, by application of both sides of the equation to $\xi$ and evaluation at the initial value $X_0$ of the microdynamics, that
\begin{equation}
x_{t+1} = \sum_{k=0}^t [P(\rho^k\circ F)](x_{t-k}) + \rho^{t+1}(X_0).
\label{eq:MZderivarion}
\end{equation}
where $\rho^k := (Q\mathcal{K})^k\xi$. The derivation of Equation \eqref{eq:MZderivarion} is explained in detail in Appendix \ref{sec:appendixMZdetails}, together with interpretation of terms of its right-hand side. 

Substituting the definition of $P$ as the orthogonal projection onto basis functions as in \eqref{eq:orthProj}, we obtain
\begin{equation}
\begin{split}
P(\rho^k\circ F)(x_{t-k}) =\varphi(x_{t-k})\langle\varphi,\varphi\rangle^{-1} \langle\varphi,\rho^k \circ F\rangle =: \varphi(x_{t-k}) h_k \in \R^m\\
\end{split}
\label{eq:MZcoeffs}
\end{equation}
with vector-valued coefficients $h_k = \langle\varphi,\varphi\rangle^{-1}\int_{\mathbb{X}}\varphi(\xi(X))^T\rho^k(F(X))d\mu(X)$.

Finding a suitable approximation of the non-accessible noise term $\rho^{t+1}(X_0)$ in \eqref{eq:MZderivarion} is generally a non-trivial task and depends on properties of the microdynamics. Examples are discussed in \cite{chu, espanol, kondrashov}. From this point onwards, we will make the simplification of replacing $\rho^{t+1}(X_0)$ by a zero-mean stochastic noise term $\varepsilon_{t+1}\in \R^m$. A typical practice is to let $\varepsilon_{t+1}$ be a zero-mean Gaussian random variable as, e.g., in~\cite{linfu, bakera}. With this, we obtain the macrodynamics
\begin{equation}
x_{t+1} = \sum_{k=0}^t \varphi(x_{t-k}) h_{k} + \varepsilon_{t+1}.
\label{eq:macrodynamics}
\end{equation}
As we can see, the evolution of the observations now depends on past terms, although the microdynamics are Markovian. For $k > 0$, the terms $[P(\rho^k\circ F)](x_{t-k})$ in Equation \eqref{eq:MZderivarion} and $\varphi(x_{t-k})$ in Equation \eqref{eq:macrodynamics} are usually referred to as \textit{memory terms}.

\subsection{Macrodynamics as a nonlinear autoregressive process}
\label{sec:macrodynamicsARprocess}
If it is reasonable to assume a sufficiently fast decay of the terms $h_{k}$ with increasing $k$, the memory terms that lie far in the past have negligible influence~\cite{horenkomolecules, shankar, chorin2000, zhu}.
In light of \eqref{eq:MZderivarion} and \eqref{eq:MZcoeffs} it is sufficient that the $\rho^k$ decay fast. To understand when this is the case, we recall $\rho^k = (Q\mathcal{K})^k\xi$, and assume the $\mathrm{range}(P) \approx \mathcal{H}$, i.e., functions parametrized by $\xi$ are well approximated by the chosen approximation space. Then, $\rho^k$ decays fast if $Q\mathcal{K}$ has a small norm, which is the case if $F$ mixes well functions that are perpendicular to $\mathcal{H}$. In other words, the dominant modes of $\mathcal{K}$ should align well with the space~$\mathcal{H}$. For quantitative statements we refer to~\cite{zhu}.

Thus, in order to obtain a feasible number of memory terms, from now one we approximate the dynamics by ending the sum in \eqref{eq:macrodynamics} with $k = p-1$ instead of $k = t$, i.e., by truncating the terms $\varphi(x_{t-p})h_p,\dots,\varphi(x_0)h_t$. Regarding the selection of an appropriate value for the \textit{memory depth} $p$ there are various methods such as Information Criteria \cite{konishi,aho} or the L-curve method \cite{lcurve}. We have thus derived a nonlinear autoregressive model (NAR) \cite{billings,huang} over $x$ given by
\begin{equation}
x_{t+1} = \sum_{k=0}^{p-1} \varphi(x_{t-k}) h_{k} + \varepsilon_{t+1}.
\label{eq:macrodynamicsNAR}
\end{equation}
with matrix-valued basis functions and vector-valued coefficients $h_k$. 

In Section~\ref{sec:sinar} we will introduce a method that identifies coefficients for NAR models in a way that is motivated by system identification methods such as Dynamic Mode Decomposition \cite{williams,dmd}, Extended Dynamic Mode Decomposition \cite{williams} or Sparse Identification of Nonlinear Dynamics \cite{sindy,sindyc}, see Fig.~\ref{fig:connections}, where the dynamics are expressed with a vector of scalar-valued basis functions and a matrix-valued coefficient. Having selected the scalar-valued basis functions $\tilde{\varphi}_1,\dots,\tilde{\varphi}_K$ and denoting $\tilde{\varphi} = [\tilde{\varphi}_1,\dots,\tilde{\varphi}_K]^T:\mathbb{Y}\rightarrow \R^K$, we thus formulate the macrodynamics
\begin{equation}
x_{t+1} = \sum_{k=0}^{p-1} H_{k} \tilde{\varphi}(x_{t-k})  + \varepsilon_{t+1},
\label{eq:macroAR}
\end{equation}
with $H_k \in \R^{m\times K}$. Although seeming like only a slight notational modification, both formulations represent different model forms. While in \eqref{eq:macrodynamicsNAR} the dynamics are expressed using different basis functions and the same coefficients across all coordinates, we will now switch to the framework in \eqref{eq:macroAR} where we select scalar-valued basis functions $\tilde{\varphi}_1,\dots,\tilde{\varphi}_L$ which are used for each coordinate while the coefficients for all coordinates can be different (the different rows of the $H_k$). In summary, for \eqref{eq:macrodynamicsNAR}, one chooses $L$ $m$-dimensional basis functions and finds $L$-dimensional coefficients while for \eqref{eq:macroAR}, one chooses $K$ $1$-dimensional basis functions and finds $(m\times K)$-dimensional coefficients.

Equation \eqref{eq:macroAR} is still consistent with the way we derive \eqref{eq:macrodynamicsNAR} through the Mori--Zwanzig formalism: Basis functions are evaluated at observations made at distinct times---no terms with mixed delays occur. In Appendix \ref{sec:appendix_2.8translation} we show how to choose basis functions and coefficients in each of the models to derive the equivalent dynamics. Please note that this does not mean that both model forms are always equivalent, as explained above. Merely, one can always choose $\tilde{\varphi}$ in dependence on $\varphi$, respectively vice versa, in a way that makes the dynamics equivalent.

\subsection{Stochastic microdynamics}
\label{sec:stochCase}

Let us consider stochastic dynamics
\begin{equation*}
X_{t+1} = F(X_t,\omega_t)
\label{eq:stochDynamic}
\end{equation*}
where $\omega_t \in \Omega$ is a random influence on $F$ which is now defined as $F:\mathbb{X}\times \Omega \rightarrow \mathbb{X}$. We will assume that the noise process $\omega_t$, $t\in\mathbb{N}$, is i.i.d.\ with law~$\mathbb{P}$. In this case we only strive to forecast the \emph{expected} macrodynamics, and define the (stochastic) Koopman operator as
\begin{equation*}
(\mathcal{K}\circ g)(X) = \mathbb{E}_{\mathbb{P}}[g(F(X,\omega))].
\end{equation*}

The spaces $\mathcal{G}$ and $\mathcal{H}$, just as the projection $P$ remain unchanged.
Naturally, to the derivation of the Mori--Zwanzig approximation we need to apply the necessary obvious modifications.
E.g., the last step in \eqref{eq:MZderivarion} now has to be modified as:
\begin{equation*}
\begin{split}
[P\mathcal{K}\rho^k](x_{t-k}) = \varphi(x_{t-k}) \langle\varphi,\varphi\rangle^{-1} \int_{\Omega}\int_X \varphi(\xi(X))^T \rho^k(F(X,\omega))d\mu(X) d\mathbb{P}(\omega).
\end{split}
\end{equation*}
We can thus obtain the identical structure of the macrodynamics as in \eqref{eq:macrodynamics} where for the computation of the coefficients $h_k$ in \eqref{eq:MZcoeffs} the expectation with respect to $\mathbb{P}$ had to be added.

\section{Sparse Identification of Nonlinear Autoregressive Models (SINAR)}
\label{sec:sinar}
We propose here a method of data-based identification for coefficients $H_k$ in \eqref{eq:macrodynamics} that is an extension of the Sparse Identification of Nonlinear Dynamics (SINDy) algorithm from \cite{sindy,sindyc,sindyclowdata}. SINDy can be used to identify the governing equations of a Markovian---in our case, discrete time---dynamical system
\begin{equation}
x_{t+1} = f(x_t) \in \R^m
\label{eq:dynsystemmarkov}
\end{equation}
from data
\begin{equation*}
\fatX = \left[
  \begin{array}{cccc}
    \vertbar &  & \vertbar \\
    x_0    & \ldots & x_{T-1}    \\
    \vertbar &      & \vertbar 
  \end{array}
\right], \fatX' =\left[
  \begin{array}{cccc}
    \vertbar &  & \vertbar \\
    x_1    & \ldots & x_{T}    \\
    \vertbar &      & \vertbar 
  \end{array}
\right], \fatX,\fatX' \in \R^{m\times T}.
\end{equation*}
We will extend this method to non-Markovian systems by applying SINDy to an extended version of $\fatX$, the \textit{Hankel} matrix
\begin{equation*}
\tilde{\fatX} = \left[
  \begin{array}{cccc}
    x_{p-1}    & \ldots & x_{T-1}    \\
    \vdots & & \vdots \\
    x_0 &\ldots& x_{T-p}
  \end{array}
\right].
\end{equation*}
In essence, this is the concept used for the Hankel-alternative view of Koopman (HAVOK) analysis from \cite{havok}, where an autoregressive model is identified on transformed coordinates obtained from a Singular Value Decomposition of the Hankel matrix from a scalar-valued observation function to separate linear from non-linear, or even chaotic, behaviour of a Markovian system. We, however, seek a formulation for the dynamics of multidimensional observations. In this section and by the choice of the name SINAR, we explicitly want to point out the connection of system identification methods for non-linear Markovian systems to their counterparts for non-linear non-Markovian systems (with finite memory these are NAR systems) that can be derived through the Mori--Zwanzig formalism from Section~\ref{sec:MZ}.
\subsection{SINDy: A short summary}
We start with a short description of SINDy~\cite{sindy}. In SINDy we try to approximate each coordinate of $f$ by a linear combination of basis functions $\theta_i:\R^{m}\rightarrow \R$ and define
\begin{equation*}
\Theta(x) = \left[
  \begin{array}{cccc}
    \theta_1(x)\\
\vdots\\
\theta_v(x)\\
  \end{array}
\right], \qquad
\Theta(\fatX) = \left[
  \begin{array}{cccc}
\theta_1(x_0) & \dots & \theta_1(x_{T-1})\\
\vdots & & \vdots\\
\theta_v(x_0) & \dots & \theta_v(x_{T-1})
\end{array}
\right], 
\end{equation*}
To this end, we fit a sparse coefficient matrix $\Xi \in \R^{m \times v}$ with rows $\Xi_i$ to the data $\fatX,\fatX'$ by solving for every row $\fatX'_i$ of $\fatX'$,
\begin{equation}
\Xi_i = \argmin\limits_{\Xi_i} \Vert \fatX_i' - \Xi_i \Theta(\fatX) \Vert_F + \lambda \Vert \Xi_i \Vert_1.
\label{eq:sindyEq}
\end{equation}
We then obtain the model
\begin{equation}
x_{t+1} \approx \Xi \Theta(x_t).
\label{eq:sindyApprox}
\end{equation}
In \eqref{eq:sindyEq} we enforce a sparsity constraint using the LASSO regression algorithm \cite{lasso} in which a regularisation term is added onto the coefficient matrix, in order to only obtain the basis functions from $\Theta$ that are dominant for the relation between $x_{t+1}$ and $\Theta(x_t)$.

The use of the 1-norm generates a sparse solution if we set $\lambda > 0 $ appropriately. Sparse models will often times be less accurate than non-sparse models. However, what we gain through a sparse right-hand side of \eqref{eq:sindyApprox} is a better interpretability of the model since only the dominant terms have been identified as influential to the dynamics. It is vital to set $\lambda$ so that the loss of accuracy is minimal compared to the gain in interpretability.

SINDy is closely related to the (first step of) the method of Dynamic Mode Decomposition (DMD) \cite{williams,dmd}, which aims at finding a linear connection between $x_t$ and $x_{t+1}$. To this end, one solves\footnote{In a second step, DMD then uses $\Xi$ from \eqref{eq:sindyEq} to uncover properties of the Koopman operator of the system. SINDy, instead, tries to explain the evolution of $x_t$ by basis functions that do not have to be linear. Still, essentially, the problem \eqref{eq:DMD} is equivalent to \eqref{eq:sindyEq} for $\Theta(x) = x$ and $\lambda = 0$.
 Further, there exists a sparse version of DMD \cite{jovanovic}, where the sparsity constraint is enforced by the additive 1-norm regularisation as in \eqref{eq:sindyEq}. Then the emerging minimization problem is the same as \eqref{eq:sindyEq} with $\Theta(x) = x$.}
\begin{equation}
 A = \argmin\limits_{A}   \Vert \fatX' - A \fatX \Vert_F.
 \label{eq:DMD}
\end{equation}

\subsection{Extending SINDy to SINAR}
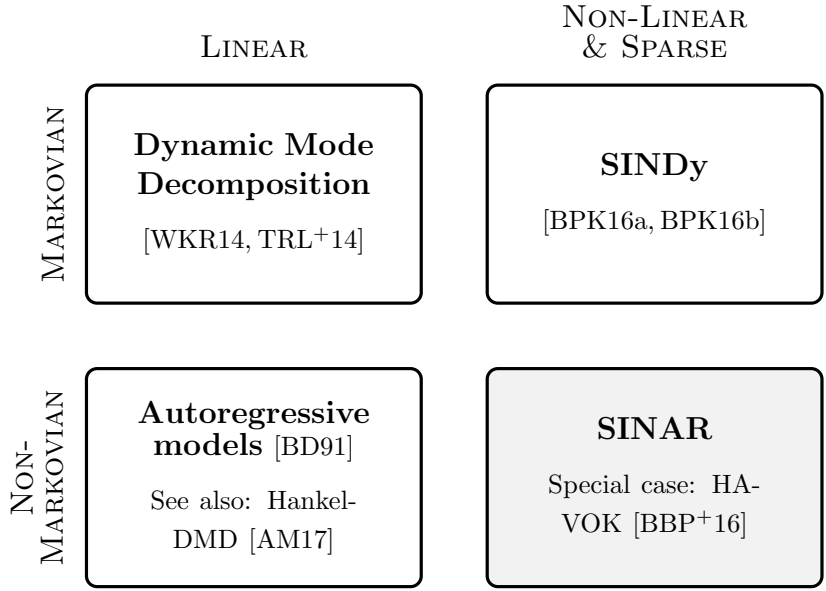
\begin{figure}[ht]
\begin{center}
\tikzstyle{block} = [rectangle, draw = black, very thick, fill=white, rounded corners, text width=10em, text centered, minimum height=7em]
	\begin{tikzpicture}
\node [block] (dmd) {\textbf{Dynamic Mode Decomposition}\\[6pt]
\footnotesize{\cite{williams,dmd}}
};
\node [block, right = 2em of dmd] (sindy) {\textbf{SINDy}\\[6pt]
\footnotesize{\cite{sindy, sindyc}}
};
\node [block, below = 2em of dmd] (ar) {\textbf{Autoregressive models} \footnotesize{\cite{brockwell}} \\[6pt]
\footnotesize{See also: Hankel-DMD \cite{hankeldmd}}
};
\node [block, right = 2em of ar, fill=gray!10] (sinar) {\textbf{SINAR}\\[6pt]
\footnotesize{Special case: HAVOK \cite{havok}}
};

\node [left = 12pt of dmd, rotate=90, anchor=center] {\textsc{{Markovian}}};
\node [left = 18pt of ar, rotate=90, anchor=center] {\stackanchor{\textsc{Non-}}{\textsc{Markovian}}};
\node [above = 6pt of dmd] {\textsc{{Linear}}};
\node [above = 6pt of sindy] {\stackanchor{\textsc{Non-Linear}}{\textsc{\& Sparse}}};

	\end{tikzpicture}
\end{center}
\caption{ Relation between different system identification methods. All of them are based on solving a Least Squares problem with respect to transformations of past to future states. While the AR minimization problem can be seen as the DMD problem on delay-embedded states and SINDy finds a nonlinear instead of linear connection between states (as in Hankel-DMD in \cite{hankeldmd}), SINAR finds a nonlinear connection between multiple past states and future ones. SINAR allows for imposing a sparsity constraint onto the determination of macro-models in the same fashion as is done in SINDy for Markovian systems. This has already been done in a special way in \cite{havok}, which is a special case of SINAR.}
\label{fig:connections}
\end{figure}

\begin{figure}[ht]
\centering
\includegraphics[width=0.9\textwidth]{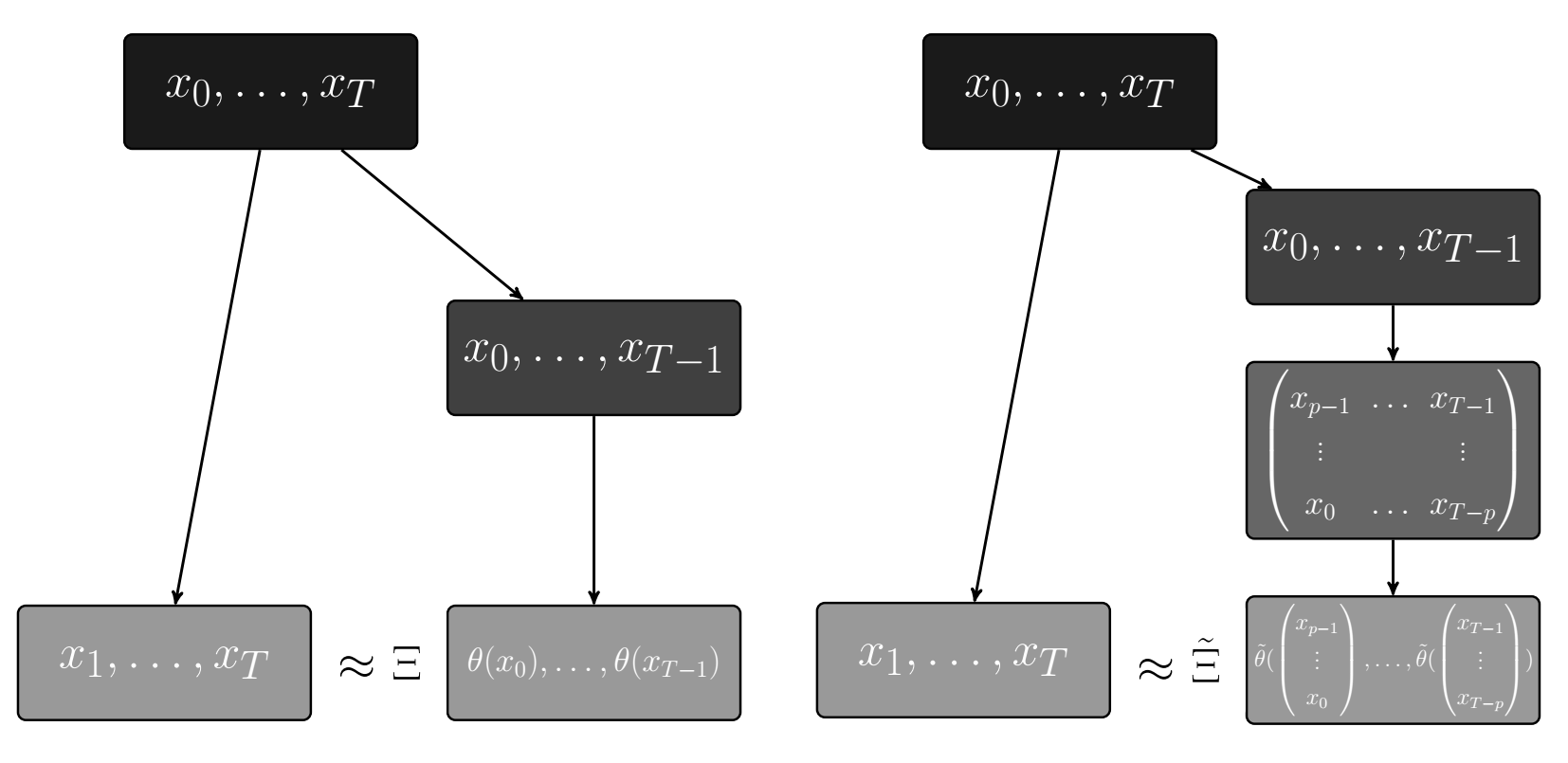}
\caption{Sketch of the SINDy algorithm (left) and SINAR (right). SINAR contains the additional step of creating a Hankel matrix.}
\label{fig:sinarbildskizze}
\end{figure}

When the dynamical model \eqref{eq:dynsystemmarkov} is insufficient in the sense that $x_{t+1}$ depends not only on $x_t$ but on memory terms too, we can apply the SINDy algorithm to suitably transformed data to obtain a nonlinear autoregressive model as in \eqref{eq:macroAR} with sparse coefficients. That is, only a few basis functions should occur with non-zero coefficients. Selecting a memory depth $p$ and denoting
\begin{equation*}
\tilde{x}_t := \left[
  \begin{array}{cccc}
    x_t\\
\vdots\\
x_{t-p+1}\\
  \end{array}
\right]
\in \R^{mp},
\end{equation*}
let us define as data matrices the Hankel matrix
\begin{equation}
\begin{split}
&\tilde{\fatX} = \left[
  \begin{array}{cccc}
    x_{p-1}    & \ldots & x_{T-1}    \\
    \vdots & & \vdots \\
    x_0 &\ldots& x_{T-p}
  \end{array}
\right] = \left[
  \begin{array}{cccc}
    \vertbar &  & \vertbar \\
    \tilde{x}_{p-1}    & \ldots & \tilde{x}_{T-1}    \\
    \vertbar &      & \vertbar 
  \end{array}
\right]\in \R^{mp\times (T-p+1)}\\
\text{ and } &\fatX' = \left[
  \begin{array}{cccc}
    \vertbar &  & \vertbar \\
    x_{p}    & \ldots & x_{T}    \\
    \vertbar &      & \vertbar 
  \end{array}
\right]\in \R^{m\times (T-p+1)}.
\end{split}
\label{eq:nardata}
\end{equation}
Again, we choose basis functions
\begin{equation*}
\tilde{\Theta}(\tilde{x}) = \left[
  \begin{array}{cccc}
    \tilde{\theta}_1(\tilde{x})\\
\vdots\\
\tilde{\theta}_v(\tilde{x})\\
  \end{array}
\right], \qquad
\tilde{\Theta}(\tilde{\fatX}) = \left[
  \begin{array}{cccc}
\tilde{\theta}_1(\tilde{x}_{p-1}) & \dots & \tilde{\theta}_1(\tilde{x}_{T-1})\\
\vdots & & \vdots\\
\tilde{\theta}_v(\tilde{x}_{p-1}) & \dots & \tilde{\theta}_v(\tilde{x}_{T-1})
\end{array}
\right]
\end{equation*}
for example 
\begin{equation*}
    \tilde{\Theta}(\tilde{x}_t) = [(x_t)_1^2,(x_t)_1 (x_t)_2,\dots, \sin((x_{t-1})_1), \dots, (x_{t-2})_m (x_{t-3})_1]^T,
\end{equation*}
and minimize for every row $\tilde{\Xi}_i$ of $\tilde{\Xi}$:
\begin{equation}
\tilde{\Xi}_i = \argmin\limits_{\tilde{\Xi}_i} \Vert \fatX'_i - \tilde{\Xi}_i \tilde{\Theta}(\tilde{\fatX}) \Vert_F + \lambda \Vert \tilde{\Xi}_i \Vert_1.
\label{eq:sinarMinProblem}
\end{equation}
Then with the basis functions with non-zero coefficients in $\tilde{\Xi}\in \R^{m\times v}$, we have derived a nonlinear autoregressive model that approximates the evolution of $x$:
\begin{equation}
x_{t+1} = \tilde{\Xi} \tilde{\Theta}(\tilde{x}_t) \in \R^{m}, \quad \text{ or, equivalently, } \quad
(x_{t+1})_i = \sum_{j=1}^{v} \tilde{\Xi}_{ij} \tilde{\theta}_j(\tilde{x}_t).
\label{eq:sindyAR}
\end{equation}
By deleting all columns of $\tilde{\Xi}$ that only contain zeros, which should be many if we enforce the sparsity constraint, we get a reduced matrix and thus a low number of terms on the right-hand side of \eqref{eq:sindyAR}. We have thus identified a sparse nonlinear autoregressive model so that we call this extension of SINDy Sparse Identification of Nonlinear Autoregressive Models (SINAR). Note that for a memory depth of $p = 1$, SINDy and SINAR are equivalent. Figure~\ref{fig:connections} shows the connections between several prominent methods for learning macro-dynamics from micro-simulation data in the Markovian and non-Markovian setting.
Figure \ref{fig:sinarbildskizze} further illustrates the different structures of SINDy and SINAR.

The choice of $\tilde{\Theta}$ allows for an arbitrary functional dependence between the distinct time-delayed observables. We can recover the special structure used in the Mori--Zwanzig formalism \eqref{eq:macrodynamicsNAR} and \eqref{eq:macroAR} by a particular choice of the basis by choosing 

\begin{equation*}
\tilde{\Theta}(\tilde{x}_t) = [
\tilde{\varphi}_1(x_t), \ldots, \tilde{\varphi}_K(x_t), \ldots, \tilde{\varphi}_1(x_{t-p+1}), \ldots, \tilde{\varphi}_K(x_{t-p+1}) ]^T,
\end{equation*}
with $\tilde{\varphi}_1,\dots,\tilde{\varphi}_K$ being scalar-valued functions as introduced in Section \ref{sec:macrodynamicsARprocess}. Then we could directly estimate the coefficients $H_k$ of the model \eqref{eq:macroAR}---which was derived through the Mori--Zwanzig formalism previously---from data, provided its distribution is approximately~$\mu$. Then $\tilde{\Xi}$ has the blockwise form
\begin{equation*}
\tilde{\Xi} = [H_0,\dots,H_{p-1}] \in \R^{m\times pK}
\end{equation*}
and
\begin{equation*}
\tilde{\Xi}\tilde{\Theta}(\tilde{x}_t) = \sum_{k=0}^{p-1}H_k \begin{bmatrix}
\tilde{\varphi}_1(x_{t-k})\\
\vdots\\
\tilde{\varphi}_K(x_{t-k})
\end{bmatrix}.
\end{equation*}
Of course, by choosing linear basis functions $\tilde{\Theta}(\tilde{x}_t) = \tilde{x}_t$ and setting $\lambda = 0$, one obtains a well-known linear autoregressive model \cite{brockwell}. Except for the sparsity term, the determination of model coefficients as in \eqref{eq:sinarMinProblem} is exactly the Least Squares method commonly used for the linear AR models. In Appendix \ref{sec:app_relationsSI}, we explain the structural equivalences and differences between SINDy, SINAR, DMD and AR models that are also sketched in Figure \ref{fig:connections}.

The covariance of the noise term $\varepsilon_{t+1}$ in \eqref{eq:macroAR} can be estimated in the common way for linear or nonlinear AR models \cite{brockwell,linfu} by calculating the statistical covariance between $\fatX'$ and $\tilde{\Xi} \tilde{\Theta}(\fatX)$ (see Appendix \ref{sec:appendixCovNAR} for more details on both statements).

In Appendix \ref{sec:henon}, we apply SINAR to an extended Hénon system, a two-dimensional dynamical system that admits a global attractor, and inspect both its accuracy in short-term predictions and its capacity to reconstruct the original attractor. This is to illustrate basic properties of non-linear autoregressive models for a simple system yielding complex dynamics.

\section{Application to an agent-based model for opinion dynamics}
\label{sec:abm}
We will now consider a network-based model of agents that change their opinions on a topic based on the opinions of their neighbors in the network. Suppose, we can only observe the percentages of agents inside the network that share each opinion, but not which agent exactly has which opinion, as in an anonymous opinion poll. Describing the evolution of these percentages can be approached by the Mori--Zwanzig formalism that we discussed in Section \ref{sec:MZ}, since they are simply observations of hidden microdynamics.
We will demonstrate the efficacy of NAR models in predicting the evolution of opinion percentages, compared with Markovian models.
We use a time-discrete agent-based model (ABM), similar to the concept of modelling opinion changes in a population explained in \cite{misra}. The ABM in \cite{misra}, however, is time-continuous while we use a time-discretized version of it. To apply the Mori--Zwanzig formalism to a time-continuous microdynamics we refer the interested reader to literature such as \cite{chorin2000,chorin2002}.

\subsection{Formulating the ABM}
The ABM is given as follows: Suppose there are $N$ agents and each agent has exactly one out of $M$ different opinions, denoted by $1,\dots,M$. The vector $X_t$, which comes from
\begin{equation*}
\mathbb{X} = \lbrace 1,\dots,M \rbrace^N,
\end{equation*}
then represents the opinions of each agent at time $t$ and $(X_t)_i$ denotes the opinion of agent $i$ at time $t$. The neighborhoods of all agents are represented by the symmetric adjacency-matrix $A \in \lbrace 0,1\rbrace^{N\times N}$ where $A_{ij} = 1$ means that agents $i$ and $j$ are neighbors of each other and $A_{ij} = 0$ otherwise. Let $N_i := \# (j : A_{ij} = 1)$ be the number of neighbors of an agent. The diagonal entries of $A$ are set to $1$, so that every agent is its own neighbour.

Let the procedure of opinion changing be given by the following rule: In every time step, every agent picks one of its neighbors in the network uniformly at random and changes its opinion with \textit{adaption probability} $\alpha_{m'm''}$ where $m'$ is the opinion of the agent and $m''$ is the opinion of the selected neighbour. This results in the term
\begin{equation*}
\mathbb{P}[(X_{t+1})_i = m'' | (X_t)_i = m'] = \alpha_{m'm''} \frac{\#( j: A_{ij} = 1 \text{ and } (X_t)_j = m'')}{N_i} \text{ for } m'\neq m'',
\end{equation*}
which we denote by $p_i^{t}(m',m'')$.
The probability for an agent not to change its opinion thus is
\begin{equation*}
p_i^t(m',m') = \mathbb{P}[(X_{t+1})_i = m' | (X_t)_i = m'] = 1-\sum_{m'' \neq m'} p_i^{t}(m',m'').
\end{equation*}
In algorithmic form, the agent-based model is executed in the following way:\\
\begin{algorithm}[H]
\SetAlgoLined
Choose end time $T$, number of agents $N$, network adjacency matrix $A$, opinion change coefficients $\alpha_{m'm''}$, initial opinions $X_0$\\
 \For{$t = 0,\dots,T$}{
 	\For{$i = 1,\dots,N$}{
 	Draw $j$ from $\lbrace j:A_{ij} = 1 \rbrace$ uniformly at random (Choose neighbour)\\
 	Draw $u_i\sim \mathcal{U}[0,1]$\\
 	If $u_i < \alpha_{(X_t)_i(X_t)_j}$: $(X_{t+1})_i = (X_t)_j$ (Adapt neighbour's opinion)
 	}
 }
 \caption{Agent-based opinion change model}
 \label{alg:ABM}
\end{algorithm}
To clarify the notation, remember that $(X_t)_i$ and $(X_t)_j$ denote the opinions of agents $i$ and $j$ at time $t$. Hence $\alpha_{(X_t)_i(X_t)_j}$ is the adaption probability of opinion $(X_t)_j$ given that an agent has opinion $(X_t)_i$. Note that in each time $t$ every agent is given the opportunity to change its opinion, and whether this happens is a probabilistic event depending only on the opinions \emph{at} time~$t$.

We can now state the so-defined microdynamics by 
\begin{equation*}
X_{t+1} = F(X_t,\omega_t)
\end{equation*}
where at every time step, $\omega_t$ denotes a tuple consisting of $N$ agents that represents the chosen neighbour of each agent plus numbers $u_i \sim \mathcal{U}[0,1]$ that govern the adaption probability $\alpha_{(X_t)_i(X_t)_j}$ as in Algorithm \ref{alg:ABM}. To be more precise, $\omega_t$ has the form
\begin{equation*}
\omega_t = [j_1,\dots,j_N,u_1,\dots,u_N],\text{ } j_i \sim \mathcal{U}\lbrace j: A_{ij} = 1 \rbrace,\text{ }  u_i \sim \mathcal{U}[0,1].
\end{equation*}
$F$ then is given by
\begin{equation*}
(X_{t+1})_i = F(X_t,\omega_t)_i = \begin{cases}
(X_t)_{j_i} & \text{ if } u_i < \alpha_{(X_t)_i,(X_t)_{j_i}} \\
(X_t)_{i} & \, \text{otherwise}.
\end{cases}
\end{equation*}
This way of stating the microdynamics seems complicated compared to the more intuitive option of denoting by $(\omega_t)_i$ the new opinion of the $i$-th agent, distributed by $[p_i^t((X_t)_i,1),\dots,p_i^t((X_t)_i,M)]$. However, this would mean that the distribution of $\omega_t$ changes over time, since the $p_i^t$ depend on $(X_t)_i$. For the Mori--Zwanzig formalism, this would prevent us from applying the procedure of skew-shift systems introduced in Section \ref{sec:stochCase} where we drew all $\omega_t$ a priori and thus independently of the $X_t$. By using the notation of $\omega_t$ denoting a tuple of neighbors $j_i$ and random numbers $u_i$ that are compared to the adaption coefficients, we can draw the whole sequence of $\omega_t$ independently of the $X_t$ and maintain consistency with the notation of skew-shift systems.

\subsection{Deducing macrodynamics from the ABM}

\paragraph{Closed-form macrodynamics.}

We now define as the \textit{opinion percentages} the function
\begin{equation*}
\xi(X)= 
\frac{1}{N}
\begin{bmatrix}
\# X_i = 1\\
\vdots\\
\# X_i = M
\end{bmatrix}
\end{equation*}
and are interested in modelling how these percentages evolve over time. It turns out that for a complete network, i.e., $A_{ij} = 1$ $\forall i,j$, we can derive macrodynamics for the expected evolution of 
\begin{equation*}
x_{t} := \xi(X_{t}),
\end{equation*}
that do not require memory terms. They are given by
\begin{equation}
\mathbb{E}[(x_{t+1})_{m'} \mid x_{t}] = (x_t)_{m'} + \sum_{m''\neq m'}(\alpha_{m''m'} - \alpha_{m'm''}) (x_t)_{m''}(x_t)_{m'} \text{ for } m' = 1,\dots,m.
\label{eq:markovMacromodel}
\end{equation}
This equation can be derived as follows: In case of a complete network, $p_i^t(m',m'') \equiv p^t(m',m'')$ is independent of $i$ because the percentages of opinions among neighbors are equal for all agents since they all have the same neighbors. Then
\begin{equation*}
p^t(m',m'') = \alpha_{m'm''} (x_t)_{m''}.
\end{equation*}
In every time step, every agent with opinion $m'$ chooses its opinion in the next time step with respective probabilities $p^t(m',m'')$ for all opinions $m'' \neq m'$ and probability $1 - \sum_{m'' \neq m'} p^t(m',m'')$ for keeping opinion $m'$. Since the number of these agents is given by $N\cdot (x_t)_{m'}$, the expected absolute number of agents that change their opinion from $m'$ to $m''$ is given by 
\begin{equation*}
\begin{split}
&\mathbb{E}[\#\text{Agents changing opinion from } m' \text{ to } m''] \\
&= \sum_{i: (X_t)_i = m'} p^t(m',m'')\\
&= N\cdot (x_t)_{m'} \cdot p^t(m',m'') \\
&= N\cdot (x_t)_{m'} \cdot \alpha_{m'm''}\cdot (x_t)_{m''}.
\end{split}
\end{equation*}
This is the expected absolute number of agents that change their opinion from $m'$ to $m''$. This means, that from this term alone, the percentage $(x_t)_{m'}$ of $m'$ is reduced by $\frac{1}{N}$ times this term, which is $\alpha_{m'm''} (x_t)_{m''} (x_t)_{m'} $. Since at the same time agents with opinion $m''$ can change their opinion to $m'$ with probability $\alpha_{m''m'} (x_t)_{m''} (x_t)_{m'} $, we have to subtract the analogous term for $\mathbb{E}[\#\text{Agents changing opinion from } m'' \text{ to } m']$ and the factor $(\alpha_{m''m'} - \alpha_{m'm''})$ comes in.
As a consequence, for a complete network the expected evolution of $x$ can be written in terms of $x$ alone, without requiring additional information of the microstate~$X$.

\paragraph{Consequences of the Mori--Zwanzig formalism.}

In the abstract language of the Mori--Zwanzig formalism from Section \ref{sec:MZ}, the above means that
\begin{equation}
    P\mathcal{K}\xi = \mathcal{K}\xi,\qquad \text{thus }Q\mathcal{K}\xi = 0,
    \label{eq:pkequalsk}
\end{equation}
because we can express $\mathcal{K}\xi = \mathbb{E}[\xi \circ F]$ as a function of $\xi$ directly by using \eqref{eq:markovMacromodel}. Let us now consider \eqref{eq:MZderivarion}, where terms of the form
\begin{equation*}
    P\mathcal{K}\rho^k \quad \text{ with } \rho^k = (Q\mathcal{K})^k\xi
\end{equation*}
occur. Equation \eqref{eq:pkequalsk} yields for $k > 0$ that $\smash{\rho^k = (Q\mathcal{K})^{k-1} (Q\mathcal{K}\xi) = 0}$.
In this way we can see that memory terms are not required for the dynamics of $\xi$ if the network is complete. However, this is generally not the case for incomplete networks, as demonstrated in detail in~\cite{banisch2014}. 
In other words, \eqref{eq:pkequalsk} is no longer valid so that the $\rho^k$ do not vanish. 
In this case, by using as $P$ the orthogonal projection onto basis functions we were able to find approximate representations of the terms $P(\rho^k\circ F)$ in~\eqref{eq:MZderivarion}. 
Here lies another part of the value of the application of the Mori--Zwanzig formalism: It installs that the structure of the ensuing macrodynamics in \eqref{eq:MZderivarion} is additive, i.e., it can be written as a sum of transformations of memory terms of \emph{individual} delays, as opposed to memory terms containing mixed delays (e.g., $\psi_1(x_t)\psi_2(x_{t-1})$). This guides our choice for a good approximation structure and reduces the number of potential basis functions from exponential in the delay depth $p$ to linear.\footnote{Supposing that there are $K$ basis functions to be used to approximate the space $\mathcal{H}$, a tensor-product basis for the complete space of ``delay-functions'' $\mathbb{Y}^p \to \mathbb{Y}$ would require $K^p$ functions. Meanwhile, the Mori--Zwanzig formalism does not mix terms from different delays, essentially working on $\bigoplus_{i=1}^p \mathcal{H}$, that is approximated by~$pK$ functions.}

For an incomplete network which is still sufficiently densely connected, we expect the microdynamics to be in expectation still close to that of a complete network. Thus, in such a case we expect $Q\mathcal{K}\xi \approx 0$, even if~\eqref{eq:pkequalsk} does not hold exactly. Consequently, assuming dense connectedness, the opinion percentages should allow for a closed-form description of their evolution with a small memory depth.
In the following, we will use SINAR to identify NAR models of this form suggested by the Mori--Zwanzig formalism.

\subsection{Recovering the macrodynamics in case of an incomplete network}
We now create realisations of the ABM with networks that consist of equally sized clusters of agents. Edges between agents from different clusters exist but are few. Inside the clusters, all agents are connected with each other. To this end, we create networks with a total number of agents $N$ consisting of equally-sized clusters. Two agents from different clusters are connected with probability $p_{between}$.

From the same initial state and with the same parameters, we create multiple realisations of the form $[X_0\dots, X_T]$ of the ABM and deduce the percentages of opinions $[x_0,\dots,x_T] = [\xi(X_0),\dots,\xi(X_T)]$.  We denote the realisations of the resulting macrodynamics by $\fatX_1,\dots,\fatX_r$ and divide these data into training data $\fatX_1,\dots,\fatX_{train}$ and validation data $\fatX_{train+1},\dots,\fatX_{r}$. Subsequently, we execute the SINAR method with different memory depths $p$ on the training data. SINAR gives us NAR models that we use for the reconstruction of the validation data. For this, the SINAR method can straightforwardly be modified for multiple trajectories by defining data matrices $\fatX' = [\fatX_1',\dots,\fatX_{train}']$ and $\tilde{\fatX} = [\tilde{\fatX}_1,\dots,\tilde{\fatX}_{train}]$ in the notation of Section \ref{sec:sinar}. We then compute the reconstruction errors of the validation data for each value of $p = 1,\dots, p_{max}$. For the reconstruction, we divide each realisation $\fatX_i$ of the validation data into blocks of length $l\geq p$. A block denotes $l$ states $\fatx^{(j)}_i= [x_{jl},\dots,x_{(j+1)l-1}]$ while the next block will be $\fatx^{(j+1)}_i= [x_{(j+1)l},\dots,x_{(j+2)l-1}].$ We then compute a reconstruction $\hat{\fatx}^{(j)}_i =[\hat{x}_{jl},\dots,\hat{x}_{(j+1)l-1}]$ of this block with the NAR model obtained with SINAR for which we use the last $p$ values of the previous block as starting values. We calculate the relative Euclidean error between reconstruction and data for each block by
\begin{equation*}
    err(\hat{\fatx}^{(j)}_i) = \frac{\Vert \fatx^{(j)}_i - \hat{\fatx}^{(j)}_i\Vert_F}{\Vert\fatx^{(j)}_i \Vert_F}.
\end{equation*}
Afterwards, we take the mean over all $err(\hat{\fatx}^{(j)}_i)$ to measure the performance of the NAR model.

Since the entries of $\xi(X_t)$ always sum up to $1$, information about the percentages of opinions $1,\dots,M-1$ immediately yields the percentage of opinion $M$ so that we use SINAR to find an NAR model for the evolution of the percentages of the first $M-1$ opinions only and omit the redundant information $\xi(X)_M$. For the reconstruction error, we compare data about the percentages of only the first $M-1$ opinions with their reconstructions. This NAR model does not necessarily ensure that the predicted first $M-1$ percentages stay between $0$ and $1$ and their sum is at most $1$. Since we make short-term predictions only, however, there will at most be only slight deviations from this property.

In the form of the diagram \eqref{eq:MZdiagram2} from Section \ref{sec:MZ}, the Mori--Zwanzig procedure applied to this concept can be described as
\begin{figure}[ht]
\centering
\includegraphics[width=0.75\textwidth]{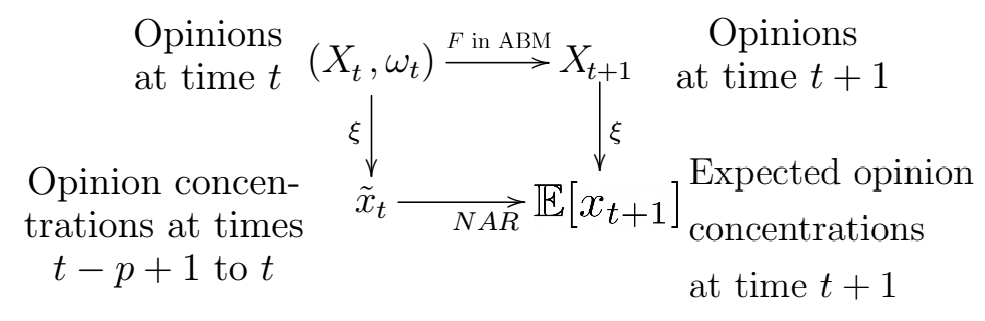}
\label{fig:diagramMZABM}
\end{figure}
\FloatBarrier
\subsubsection*{Case 1: A complete network}
For $p_{between} = 1$, the network is complete and there should be no improvement of the prediction by allowing memory terms.

We set $N = 5000, T = 300$ and $A_{ij} = 1$ $\forall i,j$. The number of different opinions is $M = 3$. As coefficients $\alpha_{m'm''}$ we choose
\begin{equation*}
\begin{bmatrix}
\alpha_{11} & \alpha_{12} & \alpha_{13}\\
\alpha_{21} & \alpha_{22} & \alpha_{23}\\
\alpha_{31} & \alpha_{32} & \alpha_{33}
\end{bmatrix}
=
\begin{bmatrix}
0 & 0.165 & 0.03\\
0.03 & 0 & 0.165\\
0.165 & 0.03 & 0
\end{bmatrix}
\end{equation*}
As initial percentages we assign values to the $(X_0)_i$ so that $\xi(X_0) = [0.45,0.1,0.45]^T$.

As the block length in the validation data, we use $l = 40$. We can already write down the macrodynamics since they are given in \eqref{eq:markovMacromodel} (see Appendix \ref{sec:appendixM3} for details):
\begin{equation}
\begin{split}
\mathbb{E}[(x_{t+1})_1 \mid x_t] &= (1+\alpha_{31} - \alpha_{13})(x_t)_1 + (\alpha_{13}-\alpha_{31}) (x_t)_1^2 + (\alpha_{21}-\alpha_{12}-\alpha_{31} + \alpha_{13}) (x_t)_1 (x_t)_2\\
&= 1.135 (x_t)_1 - 0.135 (x_t)_1^2 - 0.27 (x_t)_1 (x_t)_2,\\
\mathbb{E}[(x_{t+1})_2 \mid x_t] &= (1+\alpha_{32} - \alpha_{23})(x_t)_2 + (\alpha_{23}-\alpha_{32}) (x_t)_2^2 + (\alpha_{12}-\alpha_{21}-\alpha_{32} + \alpha_{23}) (x_t)_1 (x_t)_2\\
&= 0.865 (x_t)_2 + 0.135 (x_t)_2^2 + 0.27 (x_t)_1 (x_t)_2.
\end{split}
\label{eq:M3system}
\end{equation}
Inspired by this structure, we choose as basis functions in SINAR
\begin{equation*}
[\tilde{\varphi}_1, \dots, \tilde{\varphi}_L](x_t) = [(x_t)_1,(x_t)_2,(x_t)_1^2,(x_t)_2^2,(x_t)_1 (x_t)_2]
\end{equation*}
so that
\begin{equation}
\begin{aligned}
\tilde{\Theta}(\tilde{x}_t) &= [\underbrace{(x_t)_1,(x_t)_2,(x_t)_1^2,(x_t)_2^2,(x_t)_1 (x_t)_2}_{\text{Markovian terms as in \eqref{eq:M3system}}},\ldots \\
&\phantom{=} \underbrace{(x_{t-1})_1,(x_{t-1})_2,(x_{t-1})_1^2,(x_{t-1})_2^2,(x_{t-1})_1 (x_{t-1})_2,\dots}_{\text{Memory terms}}]^T.
\end{aligned}
\label{eq:basisfunctionsNAR}
\end{equation}
Since \eqref{eq:markovMacromodel}, resp. \eqref{eq:M3system}, describe the expected evolution of the percentages and are thus in form of deterministic models, we omit the noise term $\varepsilon_{t+1}$ from \eqref{eq:macroAR} which we assumed to satisfy $\mathbb{E}[\varepsilon_{t+1}] = 0$.

\begin{figure}[ht]
\centering
\includegraphics[width=\textwidth]{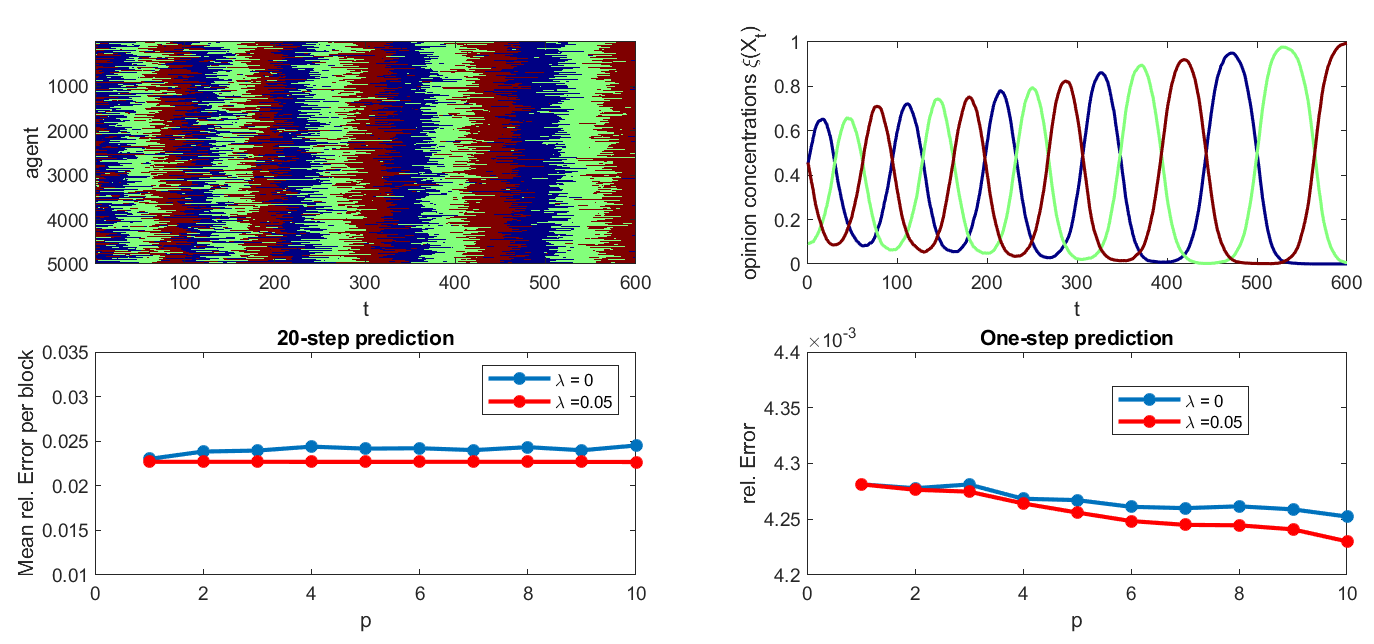}
\caption{Results for the complete network. Top left: One realisation of the microdynamics. Every column of the graphic represents the opinion of each of the $5000$ agents at one point in time. Blue denotes opinion 1, green denotes opinion 2 and red denotes opinion 3. Top right: Corresponding realisation of the macrodynamics $\xi(X)$ that represent the percentages of opinions among all agents. We can observe oscillatory behaviour since agents with opinion $1$ tend to change their opinion to $2$ and analogously from $2$ to $3$ and from $3$ to $1$. Bottom: 40-step and one-step relative prediction errors of the NAR models determined by SINAR for different memory depths $p$ with $\lambda = 0$ and $\lambda = 0.05$. As expected, the prediction error does not decrease with higher memory depth than $p=1$.}
\label{fig:N5000T300noClu1_gesamtplot}
\end{figure}

We create $r=20$ realisations of which we use $12$ for training and the others for validation. We set the sparsity parameter to $\lambda = 0$ and to $\lambda =0.05$ to test how the accuracy decreases with a sparser model. Since the macrodynamics \eqref{eq:M3system} are Markovian, we obtain for the prediction error of the validation data no improvement by allowing memory terms (Figure \ref{fig:N5000T300noClu1_gesamtplot}) for neither the 40- nor the one-step prediction error. Note that the predictions with the sparse NAR model provide slightly better accuracy for large memory depths. This is, because small non-zero coefficients for memory terms improve the fit of the training data but cause errors in the prediction of the validation data, because the macrodynamics are Markovian. Through the sparsity constraint enforced, these non-zero coefficients for memory terms are cut off. The recovered sparse macrodynamics $p = 1$ reads
\begin{equation*}
\begin{split}
(x_{t+1})_1 &= 1.1353 (x_t)_1 - 0.1351 (x_t)_1^2 - 0.2709 (x_t)_1 (x_t)_2\,,\\
(x_{t+1})_2 &= 0.8655 (x_t)_2 + 0.1344 (x_t)_2^2 + 0.2699 (x_t)_1 (x_t)_2\,,
\end{split}
\end{equation*}
which is very close to the analytically derived macrodynamics \eqref{eq:M3system}.
\FloatBarrier

\subsubsection*{Case 2: A two-cluster network}
We now construct a network with $N = 5000$ agents, divided into two clusters of size $2500$ each. We set $p_{between} = 0.0001$. Again, $M = 3$ and $\alpha_{m'm''}$ are the same as in case 1. As the starting condition, we let opinions in the first cluster be distributed by $[0.8,0.1,0.1]$ and in the second cluster by $[0.1,0.1,0.8]$. If the initial percentages in both clusters were equal then the percentages in both clusters would evolve in a quite similar way in parallel so that the macrodynamics would essentially be the same as in the complete network case. With the initial percentages being so different, it is possible that an opinion that is dominant in one cluster at one point in time but only sparsely represented in the other, can become popular through the links between agents from different clusters. This will cause the difference in behaviour of the evolution of percentages compared to the complete network.

Moreover, in order to derive the Markovian macrodynamics in Equation \eqref{eq:markovMacromodel}, we needed that the probabilities for an agent $i$ to change its opinion $(X_t)_i$ at time $t$, which we denoted by $p^t((X_t)_i,m'')$, be independent of~$i$. If the neighborhoods of different agents are generally different from each other, this is no longer the case. Especially so, if agents are distributed into different clusters, where opinion percentages might be very different. 
Thus, we cannot derive Markovian macrodynamics for this case but, in light of the Mori--Zwanzig formalism, we will need memory terms.

\begin{figure}[htb]
\centering
\includegraphics[width=\textwidth]{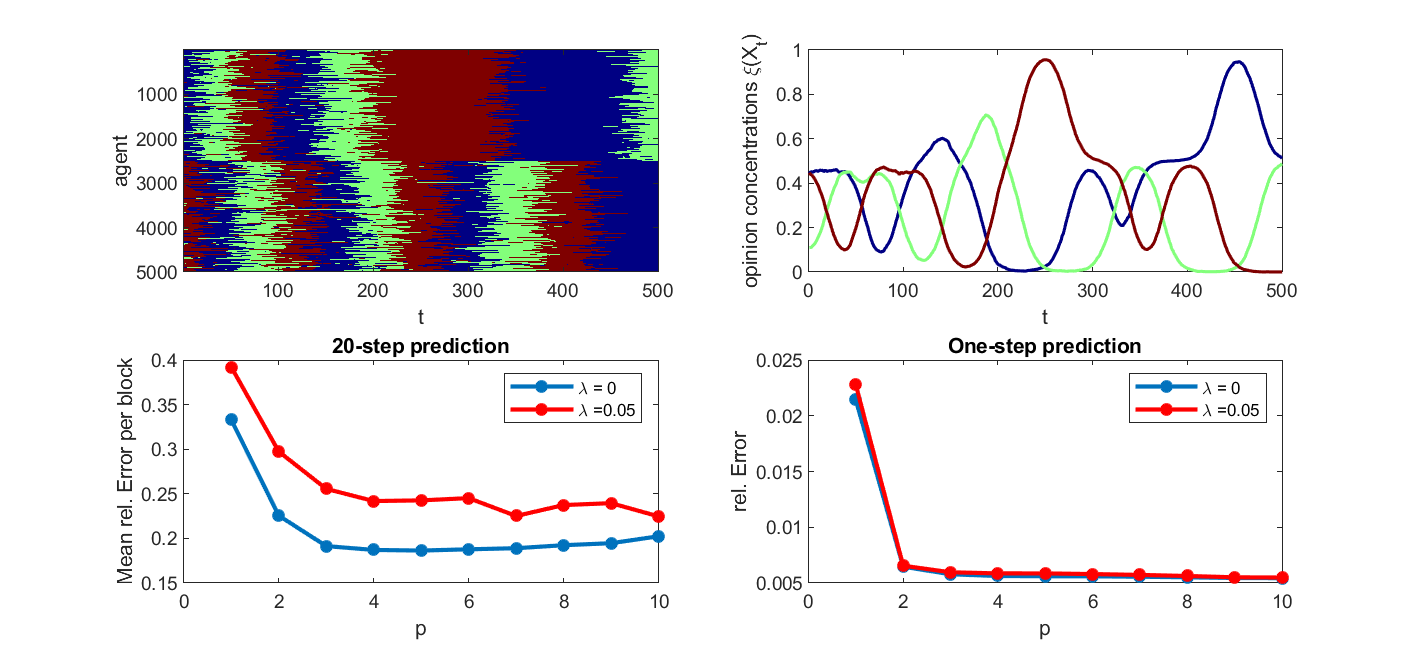}
\caption{Results for the two-cluster network. Top left: One realisation of the microdynamics. Colours represent opinions as in Figure~\ref{fig:N5000T300noClu1_gesamtplot}. Top right: Corresponding realisation of the macrodynamics $\xi(X)$. Again there is oscillatory behaviour but also plateaus and short dips as in the red and green graphs at time $25$ - $150$. This is because at these times one opinion is dominant in one cluster but not present in the other. Through the links between the clusters, an opinion, that is not present in a cluster but dominant in the other one can be revived, e.g., the blue opinion in the upper cluster. Bottom: 20-step and one-step relative prediction errors of the NAR models determined by SINAR for different memory depths $p$ with $\lambda = 0$ and $\lambda = 0.05$. Memory terms yield a significant decrease in the prediction errors compared to Markovian predictions.}
\label{fig:N5000T600noClu2_gesamtplot}
\end{figure}

To show this, we create $r = 20$ realisations of length $T = 500$ and again use $12$ for training, the remaining for validation. As block length, we choose $l = 20$. Memory terms become immediately significant, as the error graphs illustrate (Figure \ref{fig:N5000T600noClu2_gesamtplot}). We use the basis given in~\eqref{eq:basisfunctionsNAR}, which has the length~$5p$.

\begin{figure}[htb]
\centering
\includegraphics[width=\textwidth]{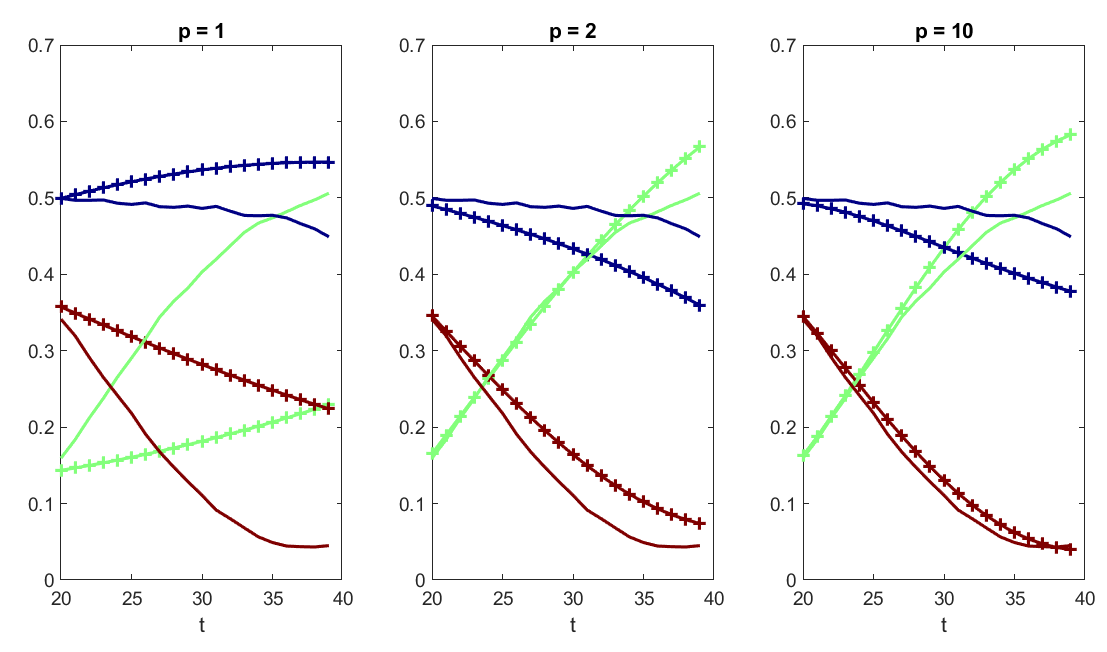}
\caption{Opinion percentages over one block of length $20$ from the validation data and prediction evolutions with NAR models obtained with SINAR for $p = 1,2$ and $10$ and $\lambda = 0$ (two-cluster network). Percentages from validation data are depicted with thin lines and predicted percentages with lines with crosses. With $p=1$, the prediction accuracy is poor and improves drastically for $p=2$. With $p=10$, the predicted evolutions come even closer to the curves from the validation data.}
\label{fig:blockPredictionnoClu2}
\end{figure}

The non-sparse and sparse solutions only deviate slightly from each other in their accuracy, but the sparse solution gives a significantly more compact model. For example, for $p = 2$, we obtain for the coefficients $\tilde{\Xi}$
\begin{equation*}
\begin{split}
\lambda = 0: \tilde{\Xi} &= \begin{bmatrix}
2.04  & 0.03 & -0.07 & -0.08 & 0.02 & -1.05 & -0.02 & 0.07 & 0.07  & -0.02 \\
-0.05 & 1.88 & 0.00  & 0.11  & 0.06 & 0.06  & -0.89 & -0.01 & -0.12 & -0.05
\end{bmatrix}\\
\lambda = 0.05: \tilde{\Xi} &= \begin{bmatrix}
1.9691 & 0      & 0 & 0 & 0 & -0.9700 & 0    & 0 & 0 & 0   \\
0      & 1.9662 & 0 & 0 & 0 & 0       & -0.9671 & 0 & 0 & 0
\end{bmatrix}
\end{split}
\end{equation*}
so that for $\lambda = 0.05$ the NAR model is given by
\begin{equation*}
\begin{split}
(x_{t+1})_1 &= 1.9691 (x_t)_1 - 0.9700 (x_{t-1})_1\\
(x_{t+1})_2 &= 1.9662 (x_t)_2 - 0.9671 (x_{t-1})_2.
\end{split}
\end{equation*}
For $p=1$, the NAR model obtained with SINAR ($\lambda = 0.05$) is
\begin{equation*}
\begin{split}
(x_{t+1})_1 &= 1.0094(x_t)_1 -0.053 (x_t)_1 (x_t)_2 \\
(x_{t+1})_2 &= 0.9894(x_t)_2 + 0.0574(x_t)_1 (x_t)_2 
\end{split}
\end{equation*}
With $\lambda = 0$, the obtained NAR model has other terms with non-zero coefficients, but these are small. In Figure \ref{fig:blockPredictionnoClu2}, an example for the predictions of opinion percentages in one block using the NAR models with $p = 1,2$ and $10$ is depicted and compared to the corresponding data. As the error graphs in Figure \ref{fig:N5000T600noClu2_gesamtplot} show already, the predicted percentages come closer to the percentages in the data with increasing memory depth.
In order to illustrate why memory terms improve the prediction accuracy, let us imagine for now that there are no links between the clusters. Then, the evolutions of opinion percentages in both clusters run in parallel to each other and are Markovian as derived previously. The opinion percentages in the full network are then given by the averages of the cluster-wise percentages $x_t^{(i)}$, i.e., $x_t = \frac{1}{2}(x_t^{(1)} + x_t^{(2)})$. This means, if we know $x_t$, then there are various options for what $x_t^{(1)}$ and $x_t^{(2)}$ can be, all of which might result in different values for $x_{t+1}^{(1)}$ and $x_{t+1}^{(2)}$ and thus $x_{t+1}$. If we are additionally given $x_{t-1}$, this might yield possible values for $x_{t-1}^{(1)}$ and $x_{t-1}^{(2)}$, which themselves make some of the candidates for $x_{t}^{(1)}$ and $x_{t}^{(2)}$ unlikely. Thus, through the information of memory terms we can restrict the options for what the percentages inside each cluster are. We illustrate this in more detail in Appendix~\ref{sec:appendixMemoryCluster}.

The links between the clusters have as consequence that within one cluster agents generally do not have identical opinion change probabilities since their neighborhoods are different. This yields additional need for memory terms since then not even for the macrodynamics in one cluster a Markovian formulation can be derived.
\subsubsection*{Case 3: A five-cluster network}
We repeat the same procedure as with the two-cluster network but with five clusters of equal size $1000$. Again, all agents within a cluster are connected with each other and $p_{between} = 0.0001$. The $\alpha_{m'm''}$ are identical to the ones used in the first two examples. As starting conditions we let opinions in the different clusters be drawn according to different distributions for each cluster. Those distributions are $[0.8,0.1,0.1],[0.1,0.1,0.8],[0.1,0.8,0.1], [0.3,0.4,0.3]$ and $[0.5,0.3,0.2]$.
The evolution of the opinion percentages is now much more irregular compared to the previous examples. The oscillatory behaviour is still present but the amplitudes differ from time to time. Through the higher number of clusters, more randomness comes into the model since an opinion can be randomly spread from one cluster, where it is dominant, to another one, where it is not dominant, suddenly altering the evolution of percentages in this cluster and thus in the whole network.

\begin{figure}[htb]
\centering
\includegraphics[width=\textwidth]{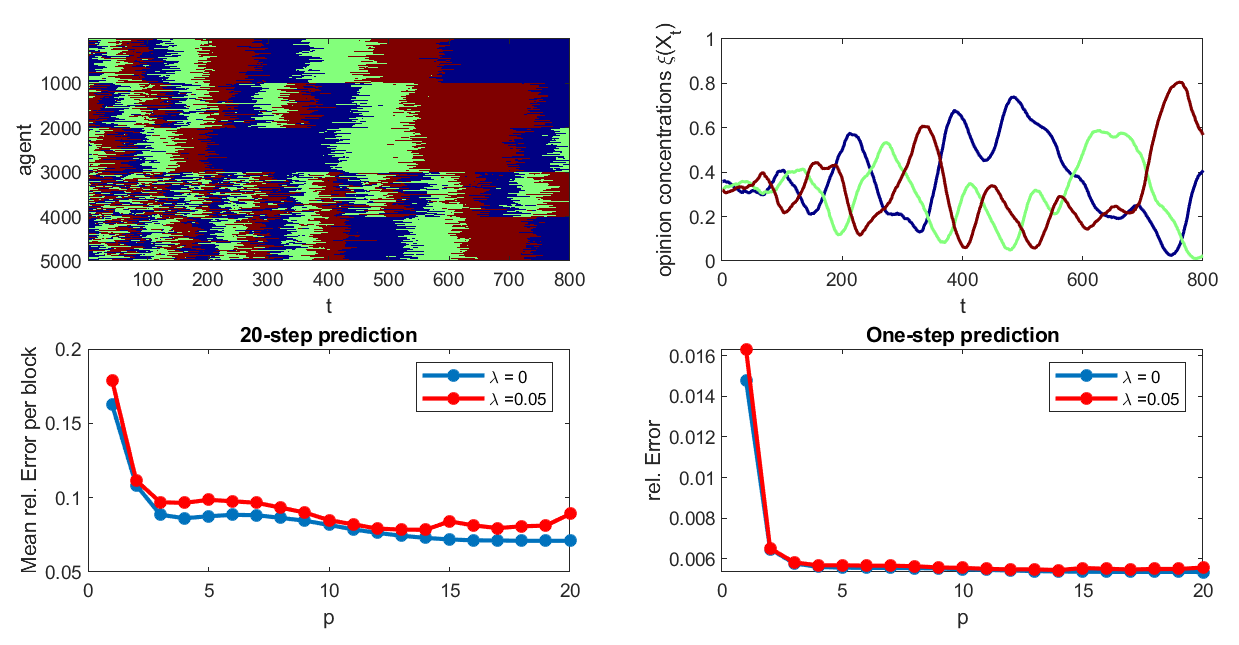}
\caption{Results for the five-cluster network. Top left: One realisation of the microdynamics. Every column of the graphic represents the opinion of each of the $5000$ agents at one point in time. Top right: Corresponding realisation of the macrodynamics $\xi(X)$. The behaviour is much more complex than in the first two cases. Bottom: 20-step and one-step relative prediction errors of the NAR models determined by SINAR for different memory depths $p$ with $\lambda = 0$ and $\lambda = 0.05$.}
\label{fig:N5000T3000noClu5_gesamtplot}
\end{figure}

We now show that, similar to when we used a two-cluster network, memory terms become important for predictions of the evolution of the microdynamics. This can be seen in Figure~\ref{fig:N5000T3000noClu5_gesamtplot}. Again, the mean relative error per block converges with increasing $p$. While in the two-cluster network example the performance did not improve visibly with $p > 10$, in this case we can get slightly lower errors for $p$ approaching $20$.

For $p = 2$ and $\lambda = 0.05$, we obtain the NAR model
\begin{equation*}
\begin{split}
(x_{t+1})_1 &= 1.8745 (x_t)_1 - 0.8748 (x_{t-1})_1\\
(x_{t+1})_2 &= 1.8672 (x_t)_2 - 0.8674 (x_{t-1})_2.
\end{split}
\end{equation*}
For $p > 2$, the models show increasing complexity, e.g., for $p=3$:
\begin{equation*}
\begin{split}
(x_{t+1})_1 &= 1.4662 (x_t)_1 - 0.1188(x_t)_2 + 0.0552 (x_t)_1^2 + 0.1318 (x_t)_1 (x_t)_2 + 0.2309 (x_{t-1})_2\\
&- 0.1899 (x_{t-1})_1(x_{t-1})_2 - 0.2021 (x_{t-1})_2^2 - 0.4658 (x_{t-2})_1 - 0.1060 (x_{t-2})_2\\
&+ 0.1206 (x_{t-2})_1^2 +0.0644 (x_{t-2})_2^2\\
(x_{t+1})_2 &= 1.3157 (x_t)_2 - 0.3161 (x_{t-2})_2.
\end{split}
\end{equation*}

\begin{figure}[ht]
\centering
\includegraphics[width=\textwidth]{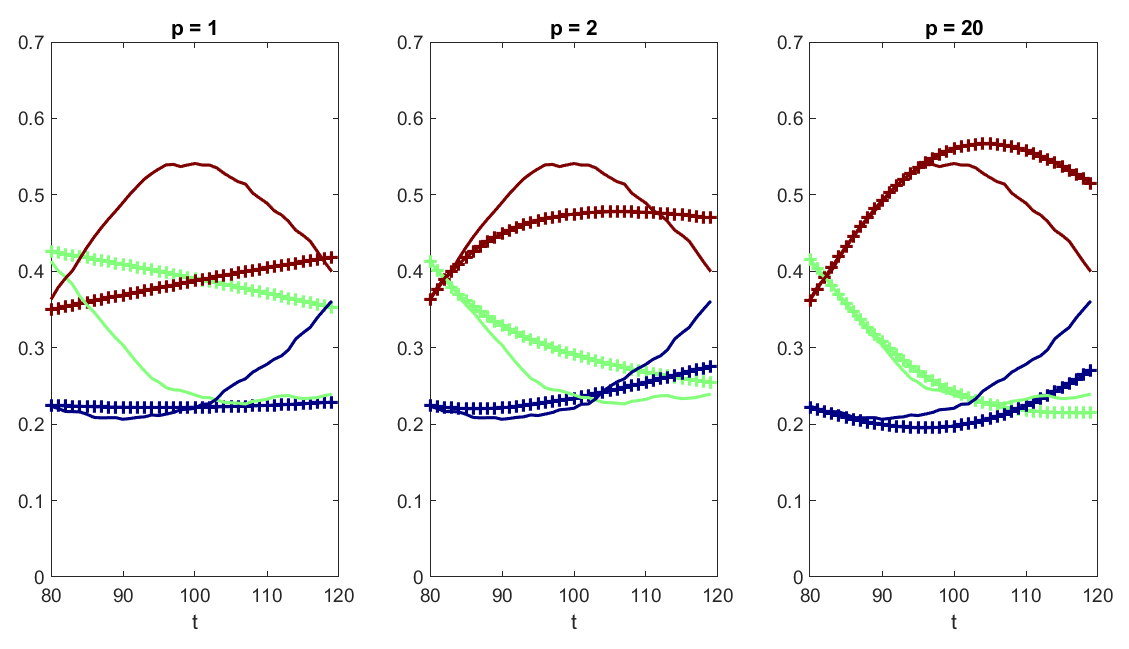}
\caption{Opinion percentages over one block of length $40$ from the validation data and prediction evolutions with NAR models obtained with SINAR for $p = 1,2$ and $20$ and $\lambda = 0$ (five-cluster network). Percentages from validation data are depicted with thin lines and predicted percentages with lines with crosses. As in the example with a two-cluster network, we can see what the error graphs in Figure \ref{fig:N5000T3000noClu5_gesamtplot} indicate: The predicted evolutions are closer to the validation data with higher memory depth of the NAR model.}
\label{fig:blockPredictionnoClu5}
\end{figure}

Again, we show as an example the predictions of percentages for one block of length 40 with memory depths $1,2$ and $10$ (Figure \ref{fig:blockPredictionnoClu5}).
As in the example with the two-cluster network, we can see that a higher memory depth indeed increases the prediction accuracy for the evolution of the opinion percentages in the short term, i.e., for predictions of length $20$ resp. $40$. Plus, enforcing the sparsity constraint with the parameter $\lambda$ in SINAR set to $0.05$ yields significantly sparser models while the prediction accuracy only suffered slightly.

\section{Discussion}
In this article, we have summarized how the evolution of observations of a dynamical system can be derived through the Mori--Zwanzig formalism, and how this can result in a nonlinear autoregressive model with memory. For the determination of model parameters, we have used methodology from data-driven system identification methods, inspired by SINDy~\cite{sindy}. We could then extend SINDy to SINAR which identifies sparse nonlinear autoregressive (NAR) models from data, thus deploying a common system identification method for non-Markovian systems.

We applied this to an agent-based model (ABM) that simulates the dynamics of opinion changes in a population. Assuming that all agents are equally strongly influenced by all other agents in the population, we showed that for the prediction of the percentages of opinions within the population memory terms are not necessary. However, for incomplete networks, this is no longer the case. Our methodology enabled us to make more accurate predictions for the percentages of opinions among the agents when the population of agents was defined by clusters with little influence between them. Additionally, sparse models obtained from enforcing a sparsity constraint in the estimation of NAR models in SINAR gave almost equally good prediction accuracy as the non-sparse ones, while yielding far simpler models. 
In the context of opinion dynamics, such sparse models permit to point out more clearly which opinions impact which others and how.

The following challenges have yet to be addressed:
\begin{itemize}
    \item In our methodology, we have assumed a noise term resulting from Mori--Zwanzig that was zero-mean. This allowed us to omit it when making predictions of the expected value of the opinion percentages. This simplifying assumption does not need to be true, and one could try to derive a more accurate representation for the noise term. As a result of this simplifiying assumption, the NAR models we considered were deterministic, even for non-deterministic microdynamics. Introduction of explicit noise in the NAR models, e.g., by extending the approach outlined in \cite{gEDMD}, could improve their (statistical) predictive capacities.
    \item  One could additionally choose a different projection $P$ in the Mori--Zwanzig formalism. The choice of an orthogonal projection on a finite set of basis functions explicitly yielded an NAR model. The right projection for a given system could inspire an optimal choice of basis functions, e.g., such that the memory depth is minimal.
    \item We have derived models that are stationary, i.e., do not change over time. Since the assumption of an equilibrium distribution over states of the microdynamics might not always hold, coefficients of the NAR model may become time-dependent. One could use a regime switching model as in \cite{horenkosocial} that fixes coefficients for a time interval before changing them to other fixed values when the macrodynamics show certain behaviour, e.g., coefficients might be different depending on which opinion is dominating.
\end{itemize}
A Matlab toolbox for the experiments done in Section \ref{sec:henon} and Appendix \ref{sec:abm} is provided under \MYhref[blue]{https://github.com/nwulkow/OpinionDyamicsModelling}{https://github.com/nwulkow/OpinionDyamicsModelling}.

\section*{Acknowledgements}

PK and CS acknowledge support by the Deutsche Forschungsgemeinschaft (DFG, German Research Foundation) under Germany's Excellence Strategy – The Berlin Mathematics
Research Center MATH+ (EXC-2046/1, project ID: 390685689). NW thanks Luzie Helfmann, Jan-Hendrik Niemann and Alexander Sikorski for helpful discussions on the subject of opinion dynamics.

\appendix
\section{Technical details on the Mori--Zwanzig equation and SINAR}

\subsection{The derivation of the Mori--Zwanzig equation}
\label{sec:appendixMZdetails}
We show here how to derive Equation \eqref{eq:MZderivarion} from the Dyson formula in Section \ref{sec:formulatingMacro}. The Dyson formula states
\begin{equation*}
\mathcal{K}^{t+1} = \sum_{k=0}^t \mathcal{K}^{t-k} P\mathcal{K} (Q\mathcal{K})^k + (Q\mathcal{K})^{k+1}.
\end{equation*}
Application of both sides of the equation to $\xi$ and evaluation at the initial value $X_0$ yields
\begin{equation}
\begin{split}
(\mathcal{K}^{t+1}\xi)(X_0) &= \sum_{k=0}^t \mathcal{K}^{t-k} [P\mathcal{K}(Q\mathcal{K})^k\xi](X_0) + (Q\mathcal{K})^{t+1}\xi(X_0)\\
\text{which results in: }  \xi(X_{t+1}) &= \sum_{k=0}^t [P\mathcal{K}(Q\mathcal{K})^k \xi](X_{t-k}) + (Q\mathcal{K})^{t+1}\xi(X_0)\\
\text{Setting } \rho^k := (Q\mathcal{K})^k\xi \text{ yields: } \xi(X_{t+1}) &= \sum_{k=0}^t [P\mathcal{K}\rho^k](X_{t-k}) + \rho^{t+1}(X_0)\\
&= \sum_{k=0}^t [P(\rho^k\circ F)](X_{t-k}) + \rho^{t+1}(X_0).
\end{split}
\label{eq:app_MZderivarion}
\end{equation}
We can replace $X_{t-k}$ by $x_{t-k}$ in the last step because the application of $P$ to a function makes this function depend only on the relevant variables. We explicitly used the parentheses around the operator $P\mathcal{K}\rho^k$ and its equivalent formulations to indicate that $P$ is a projection operator that works on the function $\mathcal{K}\rho^k$.

Since $\rho^0 = \xi$, we obtain that $P(\rho^0\circ F) = P(\xi \circ F)$. This is usually referred to as the \textit{optimal prediction} term since it is the best Markovian approximation of $\xi(X_{t+1})$, i.e., the best approximation of $\xi(X_{t+1})$ that only uses~$\xi(X_t)$. The sum in the last row of \eqref{eq:app_MZderivarion} starting at $k=1$ is referred to as the \textit{memory terms}, since these terms use information from previous values of $\xi(X)$. The term $\rho^{t+1}(X_0)$ depending on the full state $X_0$ and not on the projection $\xi(X_0)$, is often called \textit{noise}, because one does not have explicit access to it and can often only treat it as a stochastic influence.\footnote{It accumulates unobserved effects as witnessed by the complement projector~$Q$. Note that it is expected to decay fast, if the system mixes strongly (in the sense that $\mathcal{K}$ has a small spectral radius on the set of functions perpendicular to the constant function, which in turn is assumed to lie in the range of $P$). In this sense, the term 'noise' refers to negligible correlation to variables $x_{t-k}$ that contribute strongly to $\xi(X_{t+1})$.} In total, the last row of \eqref{eq:app_MZderivarion} is called the \emph{Mori--Zwanzig equation}.

Substituting the definition of $P$ as the orthogonal projection onto basis functions as in \eqref{eq:orthProj}, we obtain
\begin{equation}
\begin{split}
P(\rho^k\circ F)(x_{t-k}) &=\varphi(x_{t-k})\langle\varphi,\varphi\rangle^{-1} \langle\varphi,\rho^k \circ F\rangle\\
&= \underbrace{\varphi(x_{t-k})}_{\in \R^{m \times L}} \underbrace{\langle\varphi,\varphi\rangle^{-1}}_{\in \R^{L\times L}} \underbrace{\int_{\mathbb{X}}  \underbrace{\varphi(\xi(X))^T}_{\in \R^{L\times m}}\underbrace{\rho^k(F(X))}_{\in \mathbb{Y} \subset \R^{ m}}  d\mu(X)}_{\in \R^{L}} =: \varphi(x_{t-k}) h_k \in \R^m
\end{split}
\label{eq:appendix_MZcoeffs}
\end{equation}
with vector-valued coefficients $h_k = \langle\varphi,\varphi\rangle^{-1}\int_{\mathbb{X}}\varphi(\xi(X))^T\rho^k(F(X))d\mu(X)$.

\subsection{Translations between the model forms \eqref{eq:macrodynamicsNAR} and \eqref{eq:macroAR}}
\label{sec:appendix_2.8translation}
We show here how to translate a model in the form of \eqref{eq:macrodynamicsNAR} into the form of \eqref{eq:macroAR} and vice versa. Starting in the form of \eqref{eq:macrodynamicsNAR}, we suppose we have chosen basis functions $\varphi = [\varphi_1,\ldots,\varphi_L] \in \R^{m\times L}$ and $h_k = [(h_k)_1,\dots,(h_k)_L]\in \R^L$. This gives
\[
\varphi(x_{t-k}) h_k = \sum\limits_{i=1}^L (h_k)_i \varphi_i(x_{t-k}).
\]
Let us choose $\tilde{\varphi} = [\varphi_1^T,\ldots, \varphi_L^T]^T \in \R^{mL}$, set $H_k^{(i)} = h_i I_{m\times m} \in \R^{m\times m}$ and define $H_k = [H_k^{(1)},\dots,H_k^{(L)}]\in \R^{m\times mL}$. Then
\begin{equation*}
\begin{aligned}
    H_k \tilde{\varphi}(x_{t-k}) &= [H_k^{(1)},\dots,H_k^{(L)}]\begin{bmatrix} \varphi_1(x_{t-k})\\ \vdots\\ \varphi_L(x_{t-k}) \end{bmatrix} = \sum\limits_{i=1}^L H_k^{(i)} \varphi_i(x_{t-k}) = \sum\limits_{i=1}^L (h_k)_i \varphi_i(x_{t-k}) \\
    &= \varphi(x_{t-k}) h_k.
    \end{aligned}
\end{equation*}
Thus, we can express \eqref{eq:macrodynamicsNAR} in the form of \eqref{eq:macroAR} by imposing the restriction on the matrices $H_k$ that they have the form $H_k = h_k I_{m\times m}$. Note that we have simply modified the forms in which the dynamics are expressed but not generated a different model structure.

For the backward direction, suppose we have chosen scalar-valued basis functions $\tilde{\varphi}_1,\dots,\tilde{\varphi}_K$ and determined matrix-valued coefficients $H_k \in \R^{m\times K}$. Then we can bring \eqref{eq:macroAR} into the form of \eqref{eq:macrodynamicsNAR} by setting $L = mK$, defining $\varphi$ as the Kronecker product $\varphi = I_{m\times m} \otimes [\tilde{\varphi}_1,\dots,\tilde{\varphi}_K]$, i.e.,
\begin{equation*}
\varphi(x) = \begin{bmatrix}
\tilde{\varphi}_1(x) & \dots & \tilde{\varphi}_K(x) & 0 & & & \dots & & & 0\\
0 & \dots & 0 & \tilde{\varphi}_1(x) & \dots & \tilde{\varphi}_K(x) & & 0 & \dots & 0\\
\vdots & & & & & &\ddots & & &\\
0 & & & \dots & & & & \tilde{\varphi}_1(x) & \dots & \tilde{\varphi}_K(x)
\end{bmatrix}
\in \R^{m\times mK},
\label{eq:phiForm}
\end{equation*}
and using $mK$-dimensional coefficients
\begin{equation*}
h_k = [(H_k)_{11},\dots,(H_k)_{1K},\dots,(H_k)_{m1},\dots,(H_k)_{mK}]^T.
\end{equation*}
Then,
\begin{equation*}
\begin{aligned}
\varphi(x_{t-k}) h_k &= 
\varphi(x_{t-k})
\begin{bmatrix}
(H_k)_{11}\\
\vdots\\
(H_k)_{mK}
\end{bmatrix}
=
\begin{bmatrix}
(H_k)_{11} & \dots & (H_k)_{1K}\\
\vdots & & \vdots\\
(H_k)_{m1} & \dots & (H_k)_{mK}
\end{bmatrix}
\begin{bmatrix}
\tilde{\varphi}_1(x_{t-k})\\
\vdots\\
\tilde{\varphi}_K(x_{t-k})
\end{bmatrix} \\
&= H_k \tilde{\varphi}(x_{t-k}).
\end{aligned}
\end{equation*}

\subsection{Relation between SINDy, SINAR, DMD and AR}
\label{sec:app_relationsSI}
The diagram in Figure~\ref{fig:connections} sketches how system identification methods from different contexts are related. With DMD, SINDy, SINAR and AR models in mind, one can observe that in all of them, a minimization problem of the same form is solved: Given are data matrices $\fatX$ and $\fatX'$ which contain data points of the realisation of a (possibly memory-exhibiting) dynamical system that are shifted from each other by one time step. Then one tries to find a connection between both through a transformation of $\fatX$ which is multiplied with a coefficient matrix by solving (omitting possible sparsity constraints)
\begin{equation*}
\Xi = \argmin\limits_{\Xi} \Vert \fatX' - \Xi \Theta(\fatX) \Vert_F
\end{equation*}
In DMD, one tries to find a linear and Markovian connection between $x_t$ and $x_{t+1}$, i.e.,~$\Theta(x) =x$.
In SINDy, $\fatX$ is transformed in a possibly nonlinear way in order to explain the evolution of systems for which a linear model might be inaccurate.

Linear AR models look for a linear connection between a fixed number of past values of the system and its next value. The columns of $\fatX$, in this case, contain not just data points of the system but sequences of data points of a fixed length. More precisely, 
\begin{equation*}
    \begin{split}
        &\text{In DMD, one minimizes } \sum\limits_{t=0}^{T-1} \Vert x_{t+1} - \Xi x_t\Vert_2\\
        &\text{In AR models, define } \tilde{x}_t = [x_t^T,\dots,x_{t-p+1}^T]^T \text{and minimize } \sum\limits_{t=p-1}^{T-1} \Vert x_{t+1} - \tilde{\Xi} \tilde{x}_t \Vert_2
    \end{split}
\end{equation*}
Since in DMD one maps time-shifted versions of the same coordinates onto each other (i.e., $x_t$ to $x_{t+1}$), let us augment the AR minimazation problem to $\sum_{t=p-1}^{T-1} \Vert \tilde{x}_{t+1} - C \tilde{x}_t \Vert_2$.  Then $\tilde{\Xi}$ is equal to the upper $m$ rows of $C$ while the lower $m(p-1)$ rows of $C$ have simple structure copying the associated rows from~$\tilde{x}_t$ ($C$ is a so-called \textit{companion matrix}). In this way, the AR problem is equivalent to the DMD problem with states from the Hankel matrix defined in \eqref{eq:nardata}. In \cite{hankeldmd} the authors discuss Hankel-DMD to extract properties of the Koopman operator of a system from observational data. In doing so, they essentially fit an AR model.

In the same fashion, SINAR is the delay-embedded counterpart to SINDy and brings together SINDy and AR models in the sense that one seeks a possibly nonlinear connection between past values of the system and subsequent ones.




\subsection{Definition of the conditional expectation}
\label{sec:appendixCondExpec}
Let states $X \in \mathbb{X}$ be distributed according to $\mu$. Let us define for $\xi \in \mathcal{G}$ the level sets $L_x := \lbrace X \in \mathbb{X} : \xi(X) = x \rbrace$. Then, through the coarea formula \cite{federer}, the expectation of a function $g \in \mathcal{G}$ with $g \in L^1(\mathbb{X})$ can be written as
\begin{equation*}
\mathbb{E}_{\mu}[g(X)] = \int_{\mathbb{X}} f(X) d\mu(X) = \int_{\xi(\mathbb{X})} \int_{L_x} f(X) \mu(X) \det(\nabla \xi(X)^T \nabla \xi(X))^{-\frac{1}{2}}\, dx\, d\sigma_x(X)
\end{equation*}
where $\sigma_x$ is the Hausdorff measure on $L_x$.
Then, the conditional expectation of $f(X)$ given that $\xi(X) = x$ is (see, e.g., \cite{bittracher})
\begin{equation*}
\mathbb{E}_{\mu}[g(X)\mid \xi(X) = x] = \frac{1}{\Gamma(x)}\int_{L_x} g(X) \mu(X) \det(\nabla \xi(X)^T \nabla \xi(X))^{-\frac{1}{2}} d\sigma_x(X),
\end{equation*}
where $\Gamma(x)$ is a normalization constant.
\subsection{Determination of coefficients of linear AR models}
\label{sec:appendixLinearAR}
A linear autoregressive model with zero-mean Gaussian noise has the form
\begin{equation*}
x_{t+1} = \sum_{i = 0}^{p-1} H_i x_{t-i} + \varepsilon_{t+1},\quad  \varepsilon_{t+1} \sim \mathcal{N}(0,\Sigma^T\Sigma).
\end{equation*}
The best linear unbiased estimator (BLUE) \cite{plackett, baksalary} for the $H_i$ is the Least Squares minimizer $\tilde{\Xi} = [H_0,\dots,H_{p-1}]$, given by
\begin{equation*}
\tilde{\Xi} = \argmin_{\tilde{\Xi} = [H_0,\dots,H_{p-1}]} \Vert \fatX' - \Xi \tilde{\fatX} \Vert_F,
\end{equation*}
where $\tilde{\fatX}$ and $\fatX'$ are defined as in \eqref{eq:nardata}.

Omitting the sparsity constraint, SINAR solves the problem
\begin{equation*}
\tilde{\Xi} = \argmin_{\tilde{\Xi}} \Vert \fatX' - \tilde{\Xi} \Theta(\tilde{\fatX}) \Vert_F.
\end{equation*}
If $\Theta(\tilde{x}) = \tilde{x}$, then this is precisely the Least Squares method for linear autoregressive models.

\subsection{Covariance of noise terms of NAR models}
\label{sec:appendixCovNAR}
Assuming a relation of the form
\begin{equation*}
x'_t = \Xi \Theta(x_t) + \varepsilon_t,\quad  \varepsilon_t \sim \mathcal{N}(0,\Sigma^T \Sigma),
\end{equation*}
we find that
\begin{equation*}
Cov(x'_t - \Xi \Theta(x_t)) = Cov(\varepsilon_t).
\end{equation*}
An unbiased estimator for the covariance of a random variable $y$ is the statistical covariance
\begin{equation*}
\bar{\Sigma} = \frac{1}{T-1}\sum_{t=1}^T (y_t - \bar{y})(y_t - \bar{y})^T
\end{equation*}
where $\bar{y} = \frac{1}{T} \sum\limits_{t=1}^{T} y_t$.

In order to estimate the covariance matrix of noise terms $\varepsilon_{t+1}$ in Equation \eqref{eq:macroAR}, one has to substitute $x'_t$ by $x_{t+1}$ and $\Xi\Theta(x_t)$ by $\sum\limits_{k=0}^{p-1} H_{k} \tilde{\varphi}(x_{t-k})$ to derive the form of Equation \eqref{eq:macroAR}. Subsequently $y$ has to be substituted by $x_{t+1} - \sum\limits_{k=0}^{p-1} H_{k} \tilde{\varphi}(x_{t-k})$ and we can calculate the statistical covariance of $\varepsilon_{t+1}$ in \eqref{eq:macroAR}.

\section{Example: Application of SINAR to an extended Hénon system}

\label{sec:henon}
We demonstrate here the emergence of memory terms in the case of inaccessible variables in the sense of the Mori--Zwanzig formalism by means of an example of a dynamical system and use SINAR to detect an NAR model that reconstructs the dynamics.
\subsection{The classical Hénon system and an extension}
\begin{figure}[ht]
\centering
\includegraphics[width=0.75\textwidth]{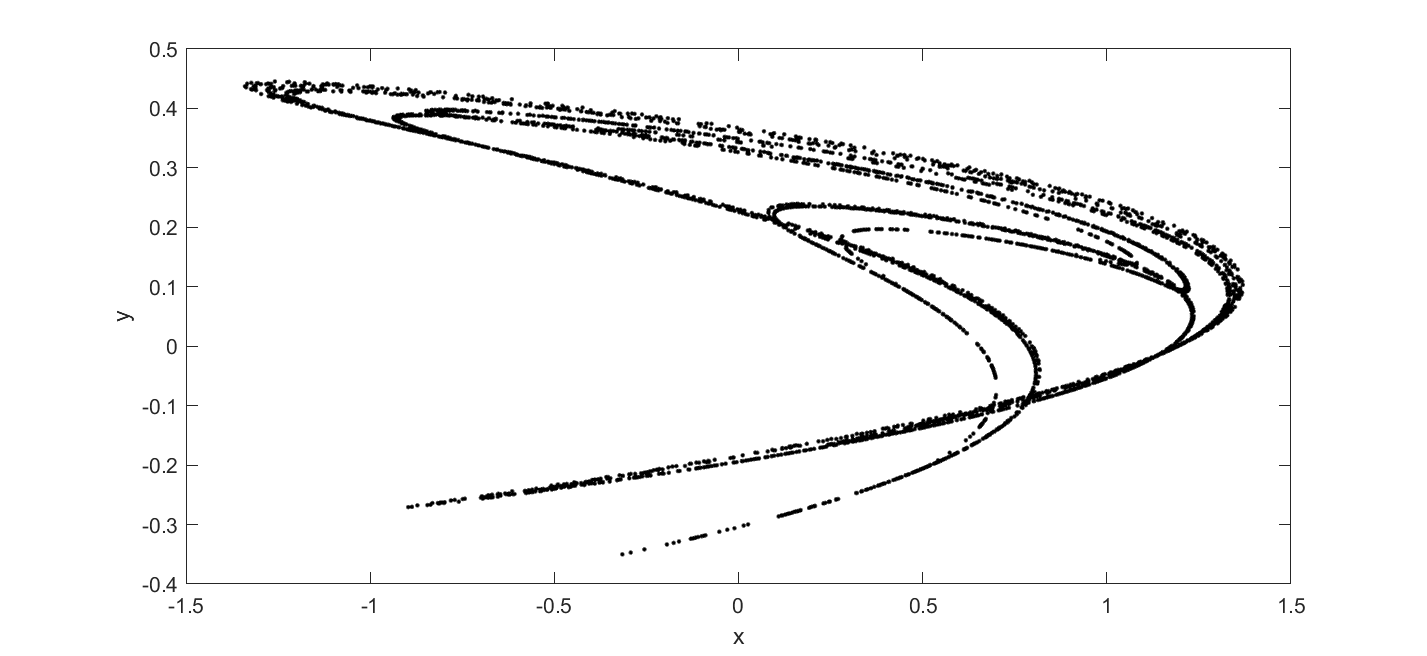}
\caption{Trajectory of length $5000$ of the two-dimensional extended Hénon system \eqref{eq:henonerweitert} with $a = 1.3, b = 0.3, c = 0.3$ and initial values $x_0 = y_0 = 0$. The first $1000$ states are omitted here so that the trajectory has time to converge towards the attractor.}
\label{fig:henonerweitertAttraktor}
\end{figure}

The classical Hénon system \cite{henon} describes a two-dimensional system that is one of the most famous examples for systems with chaotic behaviour, i.e., where slightly deviated initial conditions lead to a significantly different trajectory. The dynamical system is given by
\begin{equation*}
\begin{split}
x_{t+1} &= 1 - ax_t^2 + y_t\\
y_{t+1} &= bx_t,
\end{split}
\label{eq:henon}
\end{equation*}
where $a,b$ are fixed parameters.
As we can observe, $y_t$ is nothing more than a scaled and time-delayed version of $x_t$. We now consider $x$ as the relevant and $y$ as the irrelevant variable; this means in the Mori--Zwanzig formalism the space $\mathcal{H}$ is given by all functions depending on only $x$. We can then still express the evolution of $x$ exactly with dependence on the past two values of $x$ by plugging in the equation for $y_{t+1}$ into the equation for $x_{t+1}$:
\begin{equation*}
x_{t+1} = 1 - ax_t^2 + b x_{t-1}.
\end{equation*}
Let us now consider an extended version of the Hénon system
\begin{equation}
\begin{split}
x_{t+1} &= 1 - ax_t^2 + y_t\\
y_{t+1} &= bx_t + cy_{t}
\end{split}
\label{eq:henonerweitert}
\end{equation}
whose dynamical behaviour is visualized in Figure \ref{fig:henonerweitertAttraktor}.
\FloatBarrier
Now $y$ is more than only a scaled and time-delayed version of $x$. If we try to express $x_t$ only in dependence of its own past terms and without values of $y$ then we do not get a system with a finite memory depth but with an infinite one:
\begin{equation}
\begin{split}
x_{t+1} &= 1-ax_t^2 + bx_{t-1} + c y_{t-1}\\
&= 1- a x_{t}^2 + b x_{t-1} + cb x_{t-2} + c^2 y_{t-2}\\
&= 1- a x_{t}^2 + b x_{t-1} + cb x_{t-2} + c^2b _{t-3} +  c^3 y_{t-3}\\
&=1 - ax_{t}^2 + \sum_{j = 1}^t c^{j-1} b x_{t-j} + c^{t+1} y_0,
\end{split}
\label{eq:augmentedHenonsteps}
\end{equation}
which can be quickly shown by induction on $t$.

We have hereby derived an equation of the form of the Mori--Zwanzig equation \eqref{eq:MZderivarion} for this simple example: The term $1-ax_{t}^2$ is the optimal prediction, i.e., the Markovian approximation using the relevant variables $x_t$. The sum
\begin{equation*}
\sum_{j = 1}^t c^{j-1} b x_{t-j}
\end{equation*}
contains the memory terms depending on past values of $x$ and the term $c^t y_0$ is the noise term with information about the irrelevant, or for us inaccessible, variable $y$.
\subsection{Reconstructing the extended Hénon system with SINAR}
We now apply the SINAR algorithm to data originating from a trajectory of the extended Hénon system and demonstrate the increase in performance by using memory terms compared to applying the usual Markovian SINDy algorithm.

We set as parameters $a = 1.3, b = 0.3, c = 0.3$ and initial values $x_0 = y_0 = 0$. Then, for example, the exact model up to a memory depth of 3 in Equation \eqref{eq:augmentedHenonsteps} is
\begin{equation*}
x_{t+1} = 1 - 1.3 x_t^2 + 0.3 x_{t-1} + 0.09 x_{t-2} + 0.027x_{t-3} + \mathcal{O}(c^3).
\end{equation*}
As basis functions we choose monomials of the time-delayed coordinates up to second order without mixed terms between different delays,
\begin{equation*}
\tilde{\Theta}(\tilde{x}_t) = \left[
  \begin{array}{cccc}
1\\
x_t^2\\
x_t\\
\vdots\\
x_{t-p+1}
\end{array}
\right].
\end{equation*}

\subsubsection*{Short-term predictions}
We now generate a trajectory of length $T = 2000$ out of which we erase the first $1000$ steps to give the trajectory time to converge to the attractor. We then use the first $T_{train}$ data points for training and the remaining $1000-T_{train}$ for validation. With the training data, we determine coefficients $\tilde{\Xi}$ for the basis functions in $\tilde{\Theta}$ with SINAR for different memory depths $p$ and compute reconstructions $\hat{x}_{T_{train}+1},\dots,\hat{x}_{1000}$ of $x_{T_{train}+1},\dots,x_{1000}$ using Equation \eqref{eq:sindyAR} with initial values $x_{T_{train}-p+1},\dots,x_{T_{train}}$. In essence, we recover the coefficients of the forms $a$ resp. $c^{j-1}b$ from Equation \eqref{eq:augmentedHenonsteps} until $j = p-1$ and recompute values of the extended Hénon system with the recovered coefficients. As error measure we use the relative Euclidean prediction error
\begin{equation}
err(\hat{\fatX}') = \frac{\Vert \fatX' - \hat{\fatX}' \Vert_F}{\Vert \fatX' \Vert_F}
\label{eq:relativeError}
\end{equation}
where $\fatX' = [x_{T_{train}+1},\dots,x_{1000}]$ denotes data points from the original trajectory and $\hat{\fatX}' = [\hat{x}_{T_{train}+1},\dots,x_{1000}]$ data points from the reconstructed trajectory.

Although all coefficients are recovered up to an error of smaller than $10^{-14}$ when we use $800$ time steps for training, the reconstruction becomes inaccurate after around $100$ time steps which underlines the strongly chaotic nature of the system, i.e., small deviations at one point in time causing significant deviations in the long term behaviour. We thus use $920$ time steps for training and only $80$ time steps for validation to investigate how the relative Euclidean reconstruction error depends on the memory depth. Below we discuss how the attractor of the system is recovered using much longer reconstructions.

\begin{figure}[ht]
\centering
\includegraphics[width=1\textwidth]{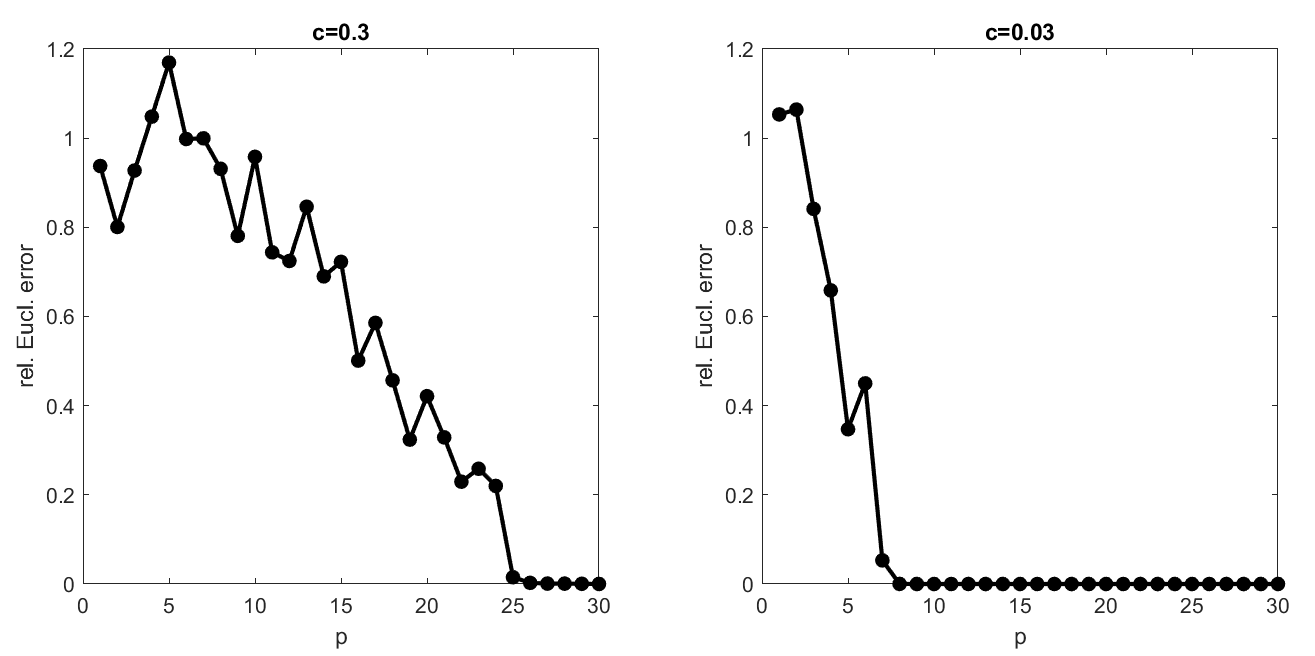}
\caption{Relative error of validation $err(\hat{\fatX}')$ for SINAR on visible variable $x$ of the extended Hénon system for two different values of $c$ with different memory depths $p$ on the x-axes. The prediction accuracy improves with increasing memory depth. Results based on SINAR with $\lambda = 0$. As parameters in the extended Hénon system, we chose $a = 1.3, b = 0.3$ and $c = 0.3$ (left) resp. $c = 0.03$ (right). For every value of $p$, the same $80$ time steps were taken into account for the reconstruction error.}
\label{fig:henonSINAR}
\end{figure}

We see in Figure \ref{fig:henonSINAR} how the relative Euclidean prediction error decreases for increasing memory depth $p$. Predicted was the evolution of $x$ with data about $x$. It is interesting to note how large a memory depth is necessary to get an accurate prediction for $x$ when $c = 0.3$ (Figure \ref{fig:henonSINAR} (left)). The chaotic nature of the system yields that even coefficients of the form $bc^{j}$ for $j = 27$ have to be taken into account. Of course, for smaller $c$ such as $c = 0.03$, memory terms in \eqref{eq:augmentedHenonsteps} decay quicker and a memory depth of $p = 8$ is sufficient to yield an accurate prediction as can be seen in Figure \ref{fig:henonSINAR} (right). For the full system $(x,y)$, the system is Markovian and the prediction error is unsurprisingly very small even for~$p = 1$.
\subsubsection*{Attractor reconstruction}
Although large deviations between original and reconstructed trajectories of $x_t$ occur after around $100$ time steps, both trajectories remain on roughly the same set of points. We quantify this by the Hausdorff distance between the two-dimensional delay embeddings  (see definition in Appendix \ref{sec:appendixHausdorff}) of the original trajectory and each reconstructed trajectory. The Hausdorff distance denotes the maximal minimal distance of members of one set of points to another set. In other words, the Hausdorff distance between two sets is $0$ if the sets are equal and big if there is a point in one set which is far away from all points in the other set.

\begin{figure}[htb]
\centering
\includegraphics[width=\textwidth]{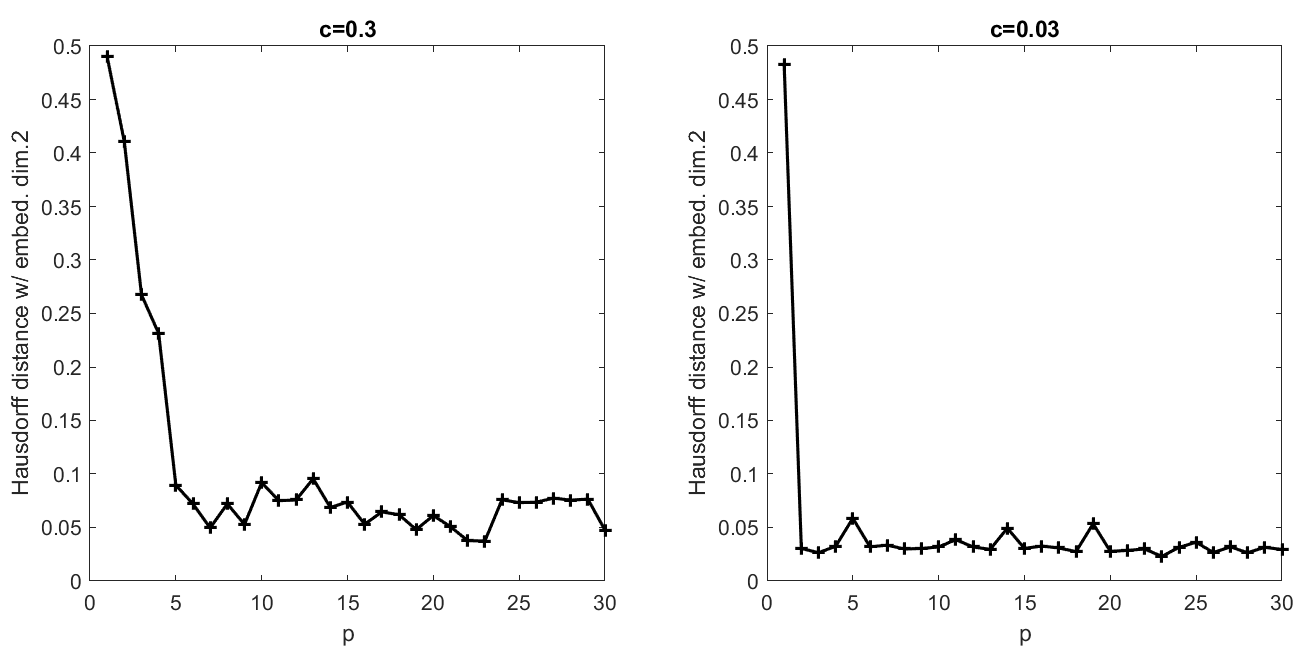}
\caption{Hausdorff distances between original and reconstructed attractors with $3000$ points of two-dimensional delay embeddings of $x$ for two different values of $c$ with different memory depths $p$ on the x-axes. Results based on SINAR with $\lambda = 0$ with parameters in the extended Hénon system $a = 1.3, b = 0.3$ and $c = 0.3$ (left) resp. $c = 0.03$ (right).}
\label{fig:henonHausdorffs}
\end{figure}

\begin{figure}[ht]
\centering
\includegraphics[width=\textwidth]{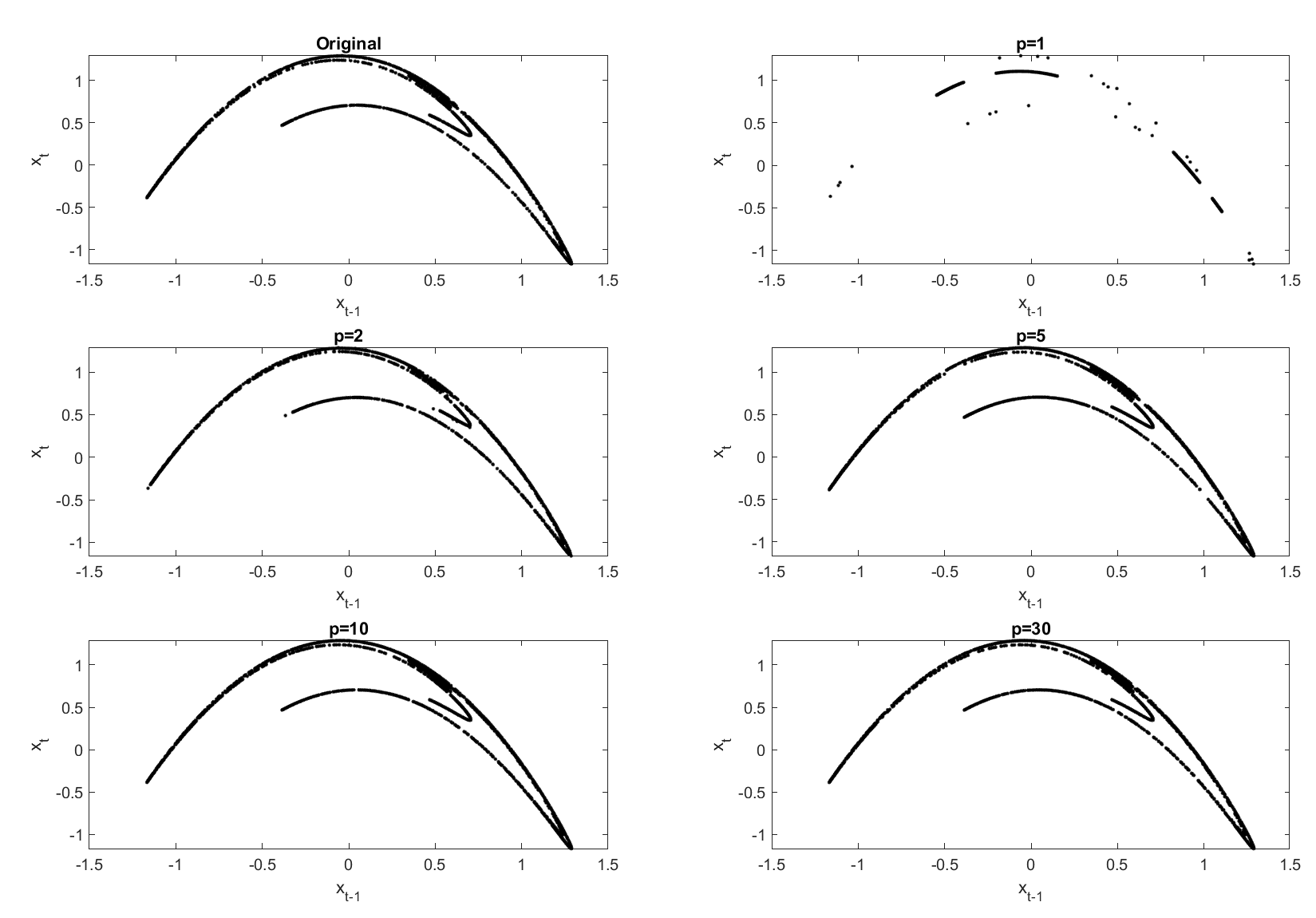}
\caption{Two-dimensional delay embedded attractors wih $3000$ points of the extended Hénon system with $a = 1.3,b = 0.3, c = 0.3$. Original (upper left) and reconstructed ones based on SINAR with $\lambda = 0$. For $p \geq 2$, differences are difficult to see but exist as the Hausdorff distances in Figure \ref{fig:henonHausdorffs} indicate.}
\label{fig:henonattractorsembed03}
\end{figure}

We make predictions of $3000$ time steps based on coefficients that were obtained with SINAR on data of $1000$ time steps. In Figure \ref{fig:henonattractorsembed03} are depicted the two-dimensional delay embeddings of the original trajectory of $x$ and the reconstructed trajectories for $p=1,2,5,10$ and $p=30$. There we see how already for $p=2$ the original and reconstructed attractors look much more similar compared to $p=1$. Figure \ref{fig:henonHausdorffs} shows the Hausdorff distances for different memory depths. Similar to the relative Euclidean prediction error, the distance decreases with increasing $p$. The remaining error is due to the fact that the complicated geometry of the attractor is hard to approximate uniformly well with a finite set of points.\footnote{Coverage of a two-dimensional object of diameter 2 by 3000 points results in a mesh size $\approx 2/\sqrt{3000} \approx 0.03$. This is the same order of magnitude as the error we observe.}

\subsection{Hausdorff distance of delay embedding of trajectories}
\label{sec:appendixHausdorff}
The Hausdorff distance between two non-empty compact sets measures the maximal minimal distance a point from one set has to the other set. It is commonly used to compare attractors of dynamical systems. The lower the Hausdorff distance between two sets, the more similar they are. 
From two trajectories $\fatX' = [x_0,\dots,x_T]$ and $\hat{\fatX}' = [\hat{x}_0,\dots,\hat{x}_T]$, we construct the delay embeddings with embedding depth $p$ as
\begin{equation*}
\mathcal{D}_p(\fatX') = \begin{bmatrix}
\begin{bmatrix}
x_{p-1}\\
\vdots\\
x_0
\end{bmatrix},
\begin{bmatrix}
x_p\\
\vdots\\
x_1
\end{bmatrix}
,\dots
\end{bmatrix},\quad \mathcal{D}_p(\hat{\fatX}') = \begin{bmatrix}
\begin{bmatrix}
\hat{x}_{p-1}\\
\vdots\\
\hat{x}_0
\end{bmatrix},
\begin{bmatrix}
\hat{x}_p\\
\vdots\\
\hat{x}_1
\end{bmatrix}
,\dots
\end{bmatrix}.
\end{equation*}
We then calculate their Hausdorff distance as
\begin{equation*}
\max( \max_{x\in \mathcal{D}_p(\fatX')} \min_{ \hat{x} \in \mathcal{D}_p(\hat{\fatX}')} \Vert x - \hat{x}\Vert_2, \max_{\hat{x}\in \mathcal{D}_p(\hat{\fatX}')} \min_{ x \in \mathcal{D}_p(\fatX')} \Vert x - \hat{x}\Vert_2 ).
\end{equation*}

\section{Details on expected opinion dynamics}
\subsection{Derivation of Equation \eqref{eq:M3system}}
\label{sec:appendixM3}
With $m=3$ opinions, Equation \eqref{eq:markovMacromodel} reads
\begin{equation*}
\begin{split}
(x_{t+1})_1 &=  (x_{t})_1 + (\alpha_{21}-\alpha_{12})(x_{t})_1 (x_{t})_2 + (\alpha_{31}-\alpha_{13})(x_{t})_1 (x_{t})_3\\
(x_{t+1})_2 &=  (x_{t})_2 + (\alpha_{12}-\alpha_{21})(x_{t})_1 (x_{t})_2 + (\alpha_{32}-\alpha_{23})(x_{t})_2 (x_{t})_3\\
(x_{t+1})_3 &=  (x_{t})_3 + (\alpha_{13}-\alpha_{31})(x_{t})_1 (x_{t})_3 + (\alpha_{23}-\alpha_{32})(x_{t})_2 (x_{t})_3.\\
\end{split}
\end{equation*}
Using $(x_{t})_3 = 1 - (x_{t})_1 - (x_{t})_2$, we get
\begin{equation*}
\begin{split}
(x_{t+1})_1 &=  (x_{t})_1 + (\alpha_{21}-\alpha_{12})(x_{t})_1 (x_{t})_2 + (\alpha_{31}-\alpha_{13})(x_{t})_1 (1-(x_{t})_1 - (x_{t})_2)\\
(x_{t+1})_2 &=  (x_{t})_2 + (\alpha_{12}-\alpha_{21})(x_{t})_1 (x_{t})_2 + (\alpha_{32}-\alpha_{23})(x_{t})_2 (1-(x_{t})_1 - (x_{t})_2).
\end{split}
\end{equation*}
Rearranging gives
\begin{equation*}
\begin{split}
(x_{t+1})_1 &= (1+\alpha_{31} - \alpha_{13})(x_t)_1 + (\alpha_{13}-\alpha_{31}) (x_t)_1^2 + (\alpha_{21}-\alpha_{12}-\alpha_{31} + \alpha_{13}) (x_t)_1 (x_t)_2\\
(x_{t+1})_2 &= (1+\alpha_{32} - \alpha_{23})(x_t)_2 + (\alpha_{23}-\alpha_{32}) (x_t)_2^2 + (\alpha_{12}-\alpha_{21}-\alpha_{32} + \alpha_{23}) (x_t)_1 (x_t)_2.
\end{split}
\end{equation*}
This is Equation \eqref{eq:M3system}.

\subsection{Representations of uncoupled expected two-cluster dynamics}
\label{sec:appendixMemoryCluster}
In this subsection, we discuss the derivation of NAR models for a network which consists of two equally-sized clusters without links between them. Having derived the expected dynamics for a complete network in Equation \eqref{eq:markovMacromodel}, we assume for now that the expected dynamics are identical with the true dynamics in order to investigate the macrodynamics if the agents behave perfectly as expected. We then get Markovian deterministic dynamics that describe the evolution of opinion percentages in each cluster. Their means are the opinion percentages in the whole network. The derivation of an NAR model for this property is analytically challenging but numerical results suggest certain structures of the macrodynamics dependent on the initial percentages.

\paragraph{Macrodynamics inside the clusters.}
Since the clusters represent complete networks of their own, we obtain for the opinion percentages $x_t^{(i)}$ inside each cluster
\begin{equation}
\begin{split}
(x_{t+1}^{(i)})_1 &= (1+\alpha_{31} - \alpha_{13})(x_t^{(i)})_1 + (\alpha_{13}-\alpha_{31}) (x_t^{(i)})_1^2 + (\alpha_{21}-\alpha_{12}-\alpha_{31} + \alpha_{13}) (x_t^{(i)})_1 (x_t^{(i)})_2\\
(x_{t+1}^{(i)})_2 &= (1+\alpha_{32} - \alpha_{23})(x_t^{(i)})_2 + (\alpha_{23}-\alpha_{32}) (x_t^{(i)})_2^2 + (\alpha_{12}-\alpha_{21}-\alpha_{32} + \alpha_{23}) (x_t^{(i)})_1 (x_t^{(i)})_2.
\end{split}
\label{eq:appendixmacro}
\end{equation}
With $x_t = \frac{1}{2}(x_t^{(1)} + x_t^{(2)})$ and denoting $a = \alpha_{31}-\alpha_{13}, b = \alpha_{21}-\alpha_{12}-\alpha_{31} + \alpha_{13},c = \alpha_{32}-\alpha_{23}, d = \alpha_{12}-\alpha_{21}-\alpha_{32} + \alpha_{23}$, this gives
\begin{equation*}
    \begin{split}
    (x_{t+1})_1 &= (1+a)\underbrace{\frac{1}{2}((x_t^{(1)})_1+(x_t^{(2)})_1)}_{(x_t)_1}-\frac{a}{2}((x_t^{(1)})_1^2 +(x_t^{(2)})_1^2)  + \frac{b}{2}((x_t^{(1)})_1 (x_t^{(1)})_2 + (x_t^{(2)})_1 (x_t^{(2)})_2)\\
    (x_{t+1})_2 &= (1+c)\underbrace{\frac{1}{2}((x_t^{(1)})_2+(x_t^{(2)})_2)}_{(x_t)_2}-\frac{c}{2}((x_t^{(1)})_2^2 +(x_t^{(2)})_2^2)  + \frac{d}{2}((x_t^{(1)})_1 (x_t^{(1)})_2 + (x_t^{(2)})_1 (x_t^{(2)})_2).
    \end{split}
\end{equation*}
Even making the simplifying assumption that $a = -c$ and $b = -d = -2a$ as is the case for the coefficients we chose for the examples, we arrive at
\begin{equation*}
    \begin{split}
    (x_{t+1})_1 &= (1+a)\underbrace{\frac{1}{2}((x_t^{(1)})_1+(x_t^{(2)})_1)}_{(x_t)_1}-\frac{a}{2}((x_t^{(1)})_1^2 +(x_t^{(2)})_1^2) - a((x_t^{(1)})_1 (x_t^{(1)})_2 + (x_t^{(2)})_1 (x_t^{(2)})_2)\\
    (x_{t+1})_2 &= (1-a)\underbrace{\frac{1}{2}((x_t^{(1)})_2+(x_t^{(2)})_2)}_{(x_t)_2}+\frac{a}{2}((x_t^{(1)})_2^2 +(x_t^{(2)})_2^2)  + a((x_t^{(1)})_1 (x_t^{(1)})_2 + (x_t^{(2)})_1 (x_t^{(2)})_2).
    \end{split}
\end{equation*}
From this, it seems impossible to find a closed Markovian expression for $x_t$. In order to understand why memory terms should help to express the evolution of $x_t$, note the following: Given $x_{t-1}$ and $x_t$, we could now find $x_{t-1}^{(1)}$ and $x_{t-1}^{(2)}$ so that these equations would yield those values for $x_t^{(1)}$ and $x_t^{(2)}$ whose average is $x_t$. This set of pairs of $x_t^{(1)}$ and $x_t^{(2)}$ would significantly be limited compared to all pairs which have this $x_t$ as their average. From these $x_t^{(i)}$, we could compute subsequent values $x_{t+1}^{(i)}$. Hence, we would have gained a more precise estimate of $x_t^{(1)}$ and $x_t^{(2)}$ and thus of $x_{t+1}$. In the stochastic ABM, the evolution of $x_t$ is originally stochastic if it represents the percentages of opinions of agents. Hence, one would not search for the $x_{t-1}^{(i)}$ that exactly yield $x_t$ but rather make this argument in terms of probabilities. We would then get different probabilities for the $x_t^{(i)}$ dependent on what $x_{t-1}$ is.

\paragraph{Simplified example: Linear dynamics inside the clusters.}

Of course, a closed expression for the evolution of $x_{t+1}$ that depends only on memory terms of $x_t$ and not on the $x_t^{(i)}$ is desirable. However, the analytical derivation of such an expression seems out of reach. Thus, as an example for much simpler macrodynamics inside each cluster, we illustrate how one can find a closed expression for the mean of two linear dynamics. For this, let
\begin{equation*}
    \begin{split}
        x_{t+1}^{(1)} &= \lambda_1 x_t^{(1)}\\
        x_{t+1}^{(2)} &= \lambda_2 x_t^{(2)}
    \end{split}
\end{equation*}
and
\begin{equation*}
    x_t = \frac{1}{2}(x_t^{(1)} + x_t^{(2)}).
\end{equation*}
Thus,
\begin{equation*}
    \begin{split}
        x_t^{(i)} &= \lambda_i^t x_0^{(i)}, \quad i=1,2\\
        \text{and } x_t &= \frac{1}{2} (\lambda_1^t x_0^{(1)}+\lambda_2^t x_0^{(2)}).
    \end{split}
\end{equation*}
Then one can observe that
\begin{equation*}
    x_{t+1} = \frac{(\lambda_1+\lambda_2)}{2}x_t - \frac{\lambda_1 \lambda_2}{2} x_{t-1}
\end{equation*}
since
\begin{equation*}
    \begin{split}
        \frac{(\lambda_1+\lambda_2)}{2}x_t - \frac{\lambda_1 \lambda_2}{2} x_{t-1} &= \frac{1}{2}[(\lambda_1+\lambda_2)(\lambda_1^t x_0^{(1)}+\lambda_2^t x_0^{(2)}) - \lambda_1\lambda_2(\lambda_1^{t-1} x_0^{(1)}+\lambda_2^{t-1} x_0^{(2)})] \\
        &= \frac{1}{2}[(\lambda_1^{t+1} x_0^{(1)}+\lambda_2^{t+1} x_0^{(2)}) + \lambda_1\lambda_2^t x_0^{(2)}+\lambda_2\lambda_1^t x_0^{(1)} - \lambda_2\lambda_1^t x_0^{(1)} - \lambda_1\lambda_2^t x_0^{(2)}]\\
        &= \frac{1}{2}(\lambda_1^{t+1} x_0^{(1)}+\lambda_2^{t+1} x_0^{(2)}) = \frac{1}{2}(x_{t+1}^{(1)} + x_{t+1}^{(2)}) = x_{t+1}. 
    \end{split}
\end{equation*}

\paragraph{Numerical results with symmetric initial percentages.}

For the macrodynamics~\eqref{eq:appendixmacro} of opinion percentages in a two-cluster network, we have not derived such a closed expression analytically. However, we can see numerically that almost exact models can be derived for a memory depth of $p=2$ if we impose \textit{symmetric} starting conditions, i.e., initial percentages that fulfill
\begin{equation*}
    (x_0^{(1)})_1 = (x_0^{(2)})_2 = 1-2(x_0^{(1)})_2, \quad (x_0^{(2)})_1 = (x_0^{(1)})_2 = 1-2(x_0^{(2)})_2.
    \label{eq:appendixStartcond}
\end{equation*}

To illustrate this, we create trajectories of length $T = 900$ of the deterministic dynamics \eqref{eq:appendixmacro} with initial percentages
\begin{equation}
    (x_0^{(1)})_1 = 0.8,\quad (x_0^{(1)})_2 = 0.1,\quad (x_0^{(2)})_1 = 0.1,\quad (x_0^{(2)})_2 = 0.8.
    \label{eq:appendixInitConc1}
\end{equation}
and $a = 0.135$ which is also the case in the examples in Section \ref{sec:abm}.

\begin{figure}[htb]
\centering
\includegraphics[width=\textwidth]{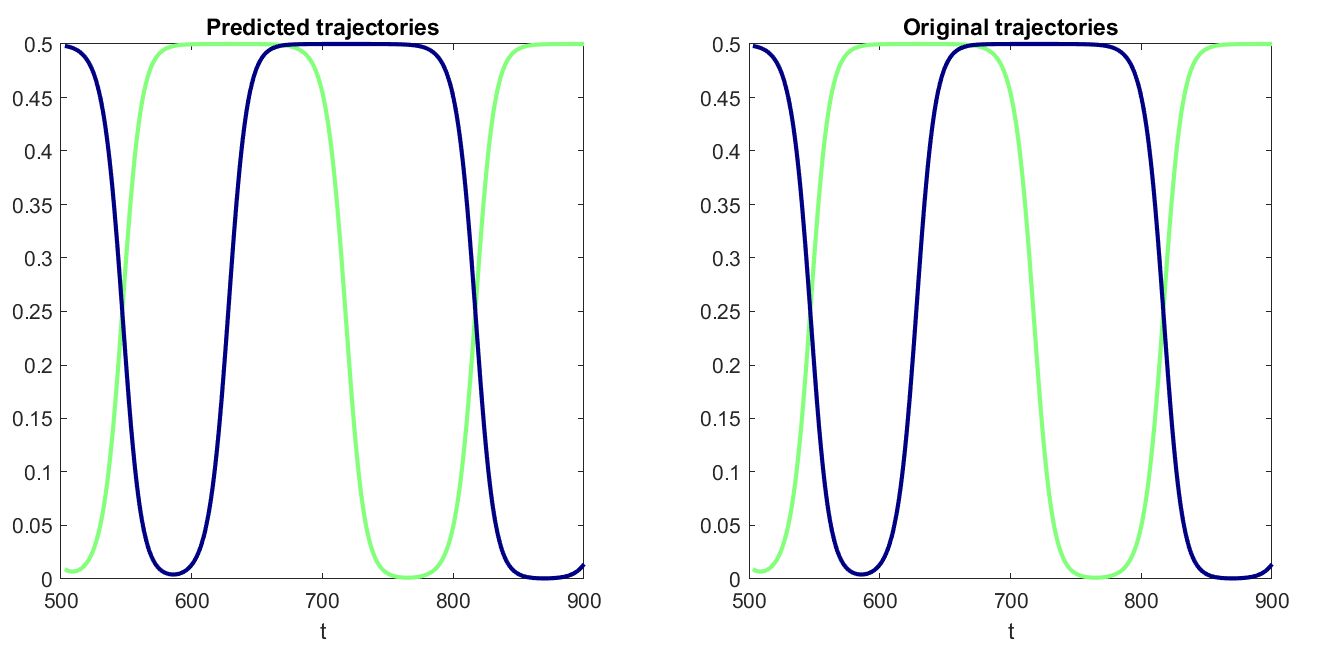}
\caption{Original trajectories for initial percentages in \eqref{eq:appendixInitConc1} and predicted trajectories with the NAR model \eqref{eq:appendixModel1}.}
\label{fig:appendix_trajectoriesModel1}
\end{figure}

From the first $500$ time steps of the resulting $x_t = \frac{1}{2}(x_t^{(1)}+x_t^{(2)})$ we estimate the NAR model (with $\lambda=0$ in SINAR)
\begin{equation}
    \begin{split}
    (x_{t+1})_1 &= 1.21 (x_t)_1 - 0.65(x_t)_2 + 0.27 (x_t)_1^2 + 0.54(x_t)_1 (x_t)_2\\
    &-0.26(x_{t-1})_1 + 0.71(x_{t-1})_2 - 0.17(x_{t-1})_1^2 - 0.10(x_{t-1})_2^2 - 0.54(x_{t-1})_1(x_{t-1})_2\\
    (x_{t+1})_2 &= -0.82 (x_t)_1 - 1.31(x_t)_2 - 0.27 (x_t)_2^2 - 0.54(x_t)_1 (x_t)_2\\
    &+0.68(x_{t-1})_1 - 0.16(x_{t-1})_2 - 0.30(x_{t-1})_1^2 - 0.03(x_{t-1})_2^2 + 0.54(x_{t-1})_1(x_{t-1})_2.
    \end{split}
    \label{eq:appendixModel1}
\end{equation}
With this model, we reconstruct the remaining $400$ time steps in the data by computing a trajectory of length $400$ with starting values given by $x_{499}$ and $x_{500}$ (Figure \ref{fig:appendix_trajectoriesModel1}). The relative Euclidean error between both trajectories amounts to $2.4\cdot 10^{-7}$. For the one-step prediction, i.e., mapping every two values $x_{t-1}$ and $x_t$ to $x_{t+1}$ with the above model, the error is $1.5\cdot 10^{-14}$. For larger memory depths, there is no improvement in prediction accuracy. This suggests that for these specific initial conditions the macrodynamics can be reproduced with memory depth $p=2$.

\paragraph{Numerical results with non-symmetric initial percentages.}

\begin{figure}[htbp]
\centering
\includegraphics[width=\textwidth]{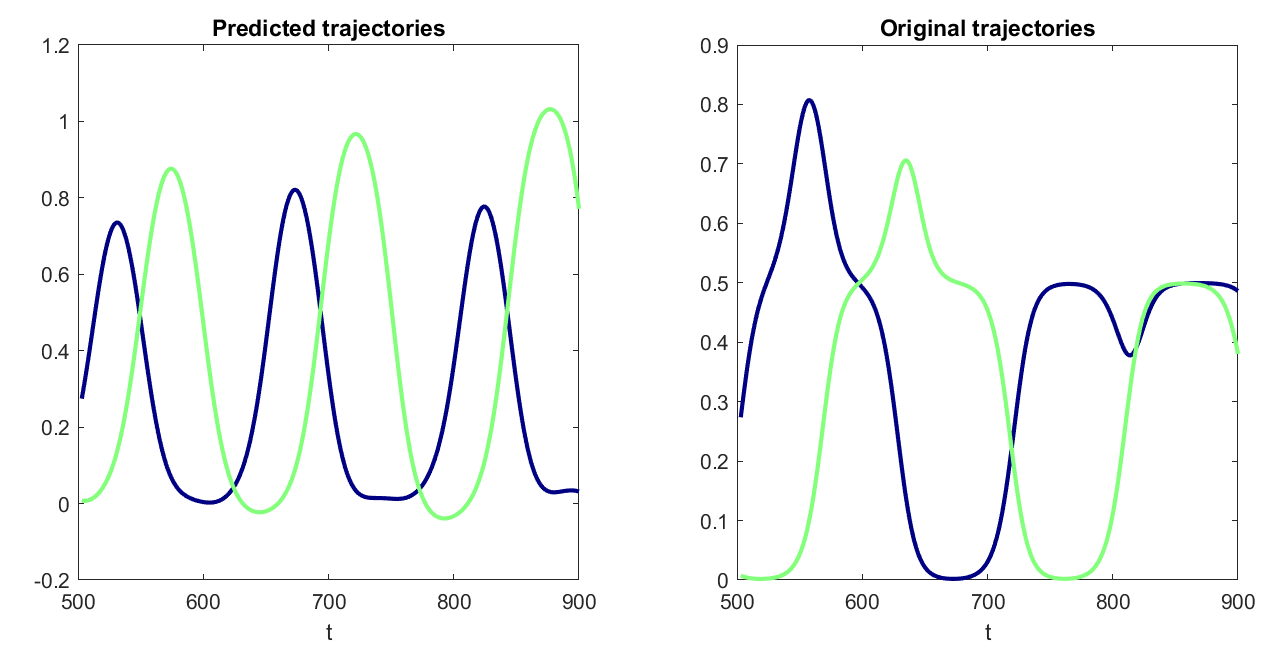}
\caption{Original trajectories for initial percentages in \eqref{eq:appendixInitConc2} and predicted trajectories with the NAR model \eqref{eq:appendixModel2}.}
\label{fig:appendix_trajectoriesModel2}
\end{figure}

For other initial percentages, we get quite different coefficients that significantly decrease the influence of the second-order terms $(x_t)_1^2,(x_t)_2^2$ and $(x_t)_1 (x_t)_2$. Let
\begin{equation}
    (x_0^{(1)})_1 = 0.7,\quad (x_0^{(1)})_2 = 0.2,\quad (x_0^{(2)})_1 = 0.1,\quad (x_0^{(2)})_2 = 0.8.
    \label{eq:appendixInitConc2}
\end{equation}
Then for $p=2$, in the same manner ($\lambda=0$), we obtain the model 
\begin{equation*}
    \begin{split}
    (x_{t+1})_1 &= 2.09 (x_t)_1 - 0.01(x_t)_2 - 0.09 (x_t)_1^2 + 0.01(x_t)_2^2- 0.15(x_t)_1 (x_t)_2\\
    &\phantom{=} -1.09(x_{t-1})_1 + 0.01(x_{t-1})_2 +0.09(x_{t-1})_1^2- 0.02(x_{t-1})_2+ 0.15(x_{t-1})_1(x_{t-1})_2\\
    (x_{t+1})_2 &= -0.04 (x_t)_1 - 1.90(x_t)_2 - 0.01 (x_t)_1^2 - 0.08(x_t)_2^2 +0.15 (x_t)_1 (x_t)_2\\
    &\phantom{=} +0.04(x_{t-1})_1 - 0.90(x_{t-1})_2 - 0.08(x_{t-1})_2^2 -0.16(x_{t-1})_1(x_{t-1})_2.
    \label{eq:appendixModel2}
    \end{split}
\end{equation*}
The original trajectories and the trajectories obtained from this model are depicted in Figure \ref{fig:appendix_trajectoriesModel2}.

\begin{figure}[htbp]
\centering
\includegraphics[width=0.9\textwidth]{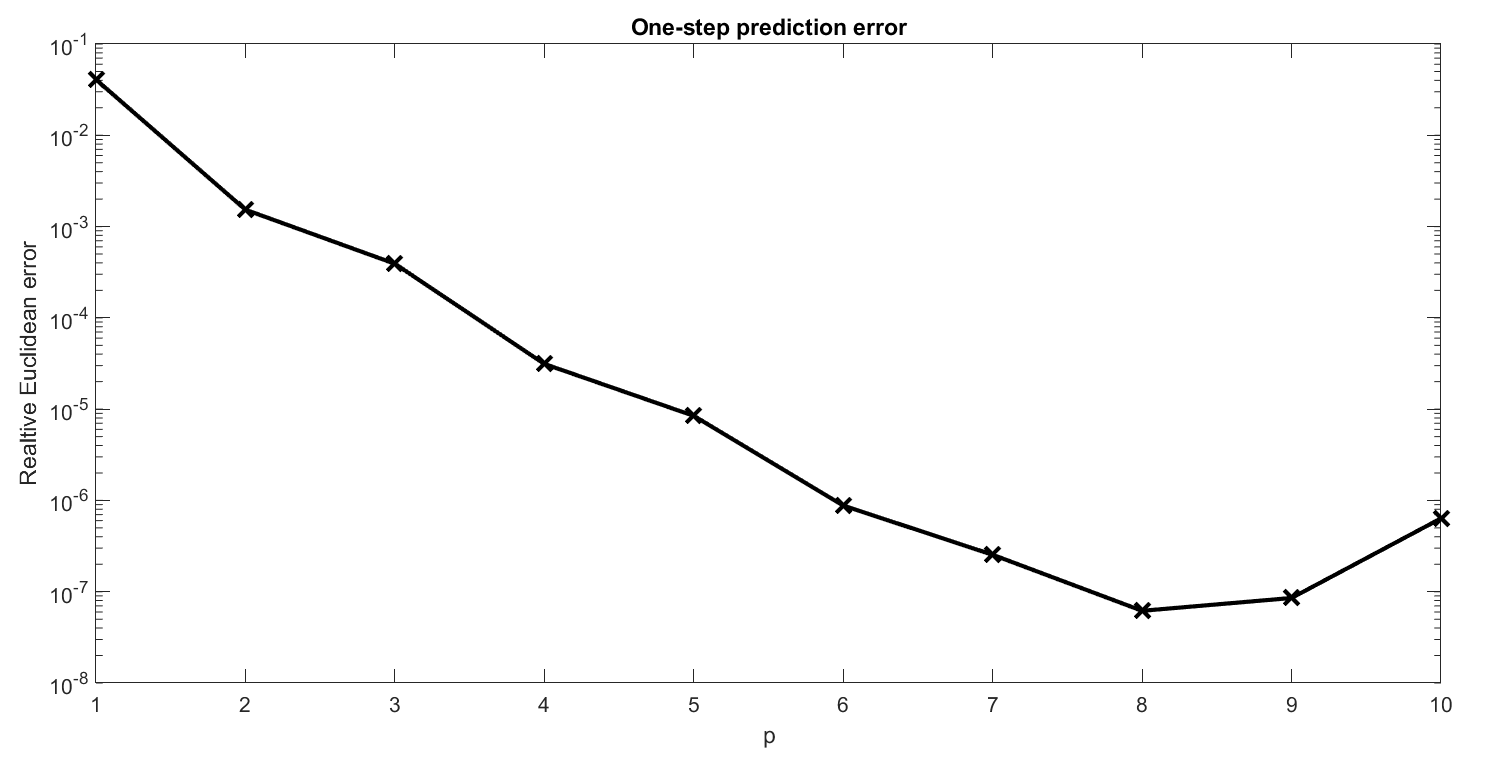}
\caption{One-step prediction error for the NAR models obtained from trajectories of $x_t$ with initial percentages of the $x_t^{(i)}$ as given in \eqref{eq:appendixInitConc2}.}
\label{fig:appendix_onesteperrorModel2}
\end{figure}

The one-step prediction error improves for memory depths larger than $p=2$ (Figure \ref{fig:appendix_onesteperrorModel2}). Since with NAR models obtained from the trajectories for these initial percentages, the predicted trajectories diverge, the full prediction error is not shown.

In summary, for a network that consists of two clusters which are uncoupled but fully connected internally, the expected macrodynamics are given by the mean of the expected intra-cluster dynamics. Assuming the dynamics to have no variance and hence to be deterministic, given in Equation \eqref{eq:appendixmacro},
with symmetric initial percentages, a memory depth of $2$ is enough for us to generate an almost exact NAR model for the macrodynamics. However, for non-symmetric initial percentages, the ensuing best-fitting NAR models with the basis functions we use are not accurate in the long-term. This seems to be in part due to the fact, that for non-symmetric initial percentages, the trajectories show more complex behaviour which no longer consists of periodic oscillations but is rather more irregular. This could cause the best-fitting NAR models to then be dominated by linear terms. 
Results about to which degree one can analytically derive NAR models for both symmetric and non-symmetric initial percentages require further research.

\small
\bibliographystyle{myalpha}
\bibliography{WulkowKoltaiSchuette_MoriZwanzigOpinionDynamics_Literature}

\end{document}